\documentclass[review]{elsarticle}

\usepackage{lineno,hyperref}
\modulolinenumbers[5]

\journal{Journal of \LaTeX\ Templates}

\usepackage{amssymb}
\usepackage{latexsym}
\usepackage{mathtools}
\usepackage{url}
\usepackage{xcolor}

\usepackage{amsmath}
\usepackage{mathtools}
\usepackage{amsfonts}
\usepackage{amsthm}
\usepackage{bm}
\usepackage{physics}

\usepackage{floatrow}

\usepackage[capitalise]{cleveref}
\graphicspath{ {images/} }

\usepackage[noend]{algpseudocode}
\usepackage{algorithm}
\usepackage[capitalise]{cleveref}
\crefformat{equation}{(#2#1#3)}
\crefname{line}{Algorithm}{algorithm}
\crefname{Item}{}{}
\usepackage{tikz-cd}

\usepackage{amsmath}
\usepackage{mathtools}
\usepackage{amsfonts}
\usepackage{amsthm}
\usepackage{bm}
\usepackage{physics}

\newtheorem{theorem}{Theorem}
\newtheorem{lemma}[theorem]{Lemma}
\newtheorem{proposition}[theorem]{Proposition}

\theoremstyle{remark}
\newtheorem{remark} {Remark}

\theoremstyle{definition}
\newtheorem{definition}{Definition}[section]

\usepackage{floatrow}

\usepackage[capitalise]{cleveref}
\graphicspath{ {images/} }

\usepackage{algpseudocode}
\usepackage{algorithm}
\usepackage[capitalise]{cleveref}
\crefformat{equation}{(#2#1#3)}
\crefname{line}{Algorithm}{algorithm}
\crefname{Item}{}{}
\usepackage{tikz-cd}

\makeatletter
\renewcommand*\env@matrix[1][*\c@MaxMatrixCols c]{%
  \hskip -\arraycolsep
  \let\@ifnextchar\new@ifnextchar
  \array{#1}}
\makeatother

\usepackage{caption}
\usepackage{subcaption}

\usepackage [autostyle, english = american]{csquotes}
\MakeOuterQuote{"}

\begin{document}

\begin{frontmatter}

\title{The \emph{curved} Mimetic Finite Difference method: allowing grids with curved faces}

\author[1]{Silvano Pitassi\corref{cor1}}
\ead{pitassi.silvano@spes.uniud.it}
\cortext[cor1]{Corresponding author:
  Tel.: +039-0432-558037;}
\author[2]{Riccardo Ghiloni}
\ead{riccardo.ghiloni@unitn.it}
\author[1]{Igor Petretti}
\ead{petretti.igor@spes.uniud.it}
\author[1]{Francesco Trevisan}
\ead{francesco.trevisan@uniud.it}
\author[1]{Ruben Specogna}
\ead{ruben.specogna@uniud.it}

\address[1]{University of Udine, Polytechnic Department of Engineering and Architecture, EMCLab, via delle scienze 206, 33100 Udine, Italy}
\address[2]{University of Trento, Mathematics Department, via Sommarive 14, 38123 Povo-Trento, Italy}


\begin{abstract}
We present a new mimetic finite difference method for diffusion problems
that converges on grids with \emph{curved} (i.e., non-planar) faces. Crucially, it gives a symmetric
discrete problem that uses only one discrete unknown per curved face. The
principle at the core of our construction is to abandon the standard definition
of local consistency of mimetic finite difference methods. Instead, we
exploit the novel and global concept of $P_{0}$-consistency.
Numerical examples confirm the consistency and the optimal convergence rate of the proposed mimetic method for cubic grids with randomly perturbed nodes as well as grids with curved boundaries.
\end{abstract}

\end{frontmatter}

\section{Introduction}
\label{sec1}

The \emph{Mimetic Finite Difference} (MFD) method
\cite{Lipnikov2014MimeticFD} is a sound numerical technique that has been
applied to solve many classes of physical problems. One of its main advantages
is the support of arbitrary \emph{polyhedral grids}, i.e., meshes composed
of general polyhedral elements with planar faces, and possibly featuring
non-matching interfaces. Polyhedral grids offer several distinctive advantages:
(i) their use can simplify the modeling of sharp geometric features of
domains; (ii) they support non-conforming mesh refinement, which does not require to trade mesh quality for accuracy;
(iii) they can be used to reduce the computational cost by mesh coarsening,
where multiple elements from a (fixed) background mesh are coalesced into
a single one.

However, many applications require grids whose elements have
\emph{curved} (i.e., non-planar) faces.
For example, in the modeling of complex reservoir geometries in real-world geological applications or of geometries like cylinders and spheres which are of common use in electromagnetics.
We emphasize that even if
the computational domain is a polyhedron, curved faces may appear in the
interior of the domain just because a hexahedral mesh has been obtained
with an unstructured meshing algorithm.

Unfortunately, the standard MFD method does not converge on grids with
curved faces \cite{Lipnikov2006TheMF}. In order to deal with such convergence
problems, the different numerical methods proposed in literature implement
two basic strategies, which, however, exhibit important limitations and
downsides.

The first approach is to approximate curved faces with triangles (or more
generally with polygons) to obtain a polyhedral grid where all the elements
have planar faces so that we can apply the standard MFD method
\cite{Lipnikov2006TheMF,brezzi2005family}. With this approach we obtain
a symmetric discrete problem but the price to pay is the additional number
of discrete face unknowns, which will be proportional to the number of
triangles partitioning every curved face.

The second approach proceeds in two steps: (i) a constant vector field
is obtained by interpolating discrete face unknowns on every element by
using standard least squares methods; (ii) mimetic inner products are constructed
on every element by projecting these constant vector fields onto edges
of a \emph{dual} (or \emph{secondary}) grid $\tilde{K}$
\cite{pitassi2021explicit}. A dual grid $\tilde{K}$ is constructed by duality
starting from a given discretization grid $K$, hence, called a
\emph{primal} grid. Here, duality is expressed as a bijective correspondence
between geometric elements of the pair of grids $(K,\tilde{K})$ such that
to each $d$-dimensional geometric element in $K$ corresponds a unique dual
$(3-d)$-dimensional geometric element in $\tilde{K}$. However, using this
dual grid approach on curved meshes the resulting discrete problem is non-symmetric
which significantly reduces the number of available efficient solution
methods. Basically, the discretization methods in
\cite{Mishev2002NonconformingFV,Aavatsmark2002AnIT} reduce to the above
approach although a dual grid is not explicitly introduced.

The aim of this paper is to introduce the \emph{curved} MFD method, an extension
of the MFD method that converges on grids with curved faces. Crucially,
it gives a discrete problem that is both symmetric and uses only one discrete
face unknown for each curved face. Our work thus answers to an open question
raised in \cite{Brezzi2006CONVERGENCEOM,Lipnikov2006TheMF} on ``whether
the use of additional degrees of freedom is the only way to preserve symmetry
in the discrete problem''. To the best of our knowledge, the numerical
method presented in \cite{Brezzi2006CONVERGENCEOM,brezzi2007new} is the
only low-order mimetic numerical method that converges on grids with curved
faces. In particular, it results in a symmetric discrete problem which
uses three discrete unknowns for each curved face (more precisely,
\emph{strongly} curved faces\footnote{\cite{Brezzi2006CONVERGENCEOM} distinguishes
between \emph{moderately} and \emph{strongly} curved faces; the distinction
between the two types of curved faces is based, as the names suggest, on
a measure of the curvature of faces and is used as a definition for the theoretical
analysis of convergence.}). More recently, high-order numerical methods
have been developed for the treatment of 3-dimensional grids with curved
faces in the context of Hybrid High Order methods \cite{yemm2023new},
\cite{Botti2018AssessmentOH} and Virtual Element Method
\cite{dassi2022bend}, but we note that their lowest order versions employ
more than one discrete unknown per curved face and thus are not equivalent
to the curved MFD method proposed in this paper. Finally, the mimetic method
developed in \cite{pitassi2021explicit} uses only one discrete unknown
per curved face but is restricted to grids obtained by barycentric subdivision
of tetrahedral grids, which are made of barycentric dual cells with non-planar
faces.

The core idea of our method is to abandon the definition of
\emph{consistency} used in the standard MFD method. Consistency is a local
property of \emph{mimetic inner products} and
\emph{reconstruction operators} that encodes an exactness property for constant
vector fields \cite{Lipnikov2014MimeticFD}. In particular, the consistency
of the standard MFD method is equivalent to the consistency enforced on
the \emph{barycentric dual grid}. In fact, as demonstrated in
\cite{Pitassi2021TheRO}, and also as we will show in
Section~\ref{comparison}, reconstruction operators of the standard MFD method
are equivalent to reconstruction operators of the Discrete Geometric
Approach (DGA) \cite{cst}. The latter, are defined by geometric elements of the barycentric dual grid and satisfy two key properties \cite{Pitassi2021TheRO}:
\begin{description}
\item[(P1)] They provide a left inverse for local projection operators
restricted to the vector subspace of constant vector fields, which is the
\emph{accuracy} or \emph{unisolvence} property for constant vector fields
\cite{da2014mimetic}.
\item[(P2)] Their geometric entries form closed paths supported on different
elements.
\end{description}
The key observation is that the property \textbf{(P1)} is not valid for
grids with curved faces: this is the theoretical reason for the inadequacy
of the consistency of the standard MFD method on grids with curved faces.

Instead of focusing on local consistency, in our curved MFD method we employ
the novel concept of $P_{0}$-consistency introduced in
\cite{Pitassi2021TheRO}. In contrast with consistency of the standard MFD
method, $P_{0}$-consistency is a global property as meaning that involves
reconstruction operators defined on more than one
grid element. The idea
underlying its definition is to abstract the properties \textbf{(P1)} and
\textbf{(P2)} of the barycentric dual grid and reverse our line of reasoning:
any grid whose geometric elements satisfy properties \textbf{(P1)} and
\textbf{(P2)} leads to consistent reconstruction operators, precisely,
$P_{0}$-consistent reconstruction operators. We will present their formal
definition in Section~\ref{math}. As a side note, we would like to point that
special cases of $P_{0}$-consistent reconstruction operators have already
been successfully employed in the mimetic scheme in
\cite{pitassi2021explicit}.

In Section~\ref{curved}, we provide a geometric characterization of $P_{0}$-consistent
reconstruction operators as the affine solution space of a linear system
of equations. Then, following \cite{brezzi2007new}, we focus on the practically
important class of cubic grids with randomly perturbed nodes (according
to a suitable probability distribution), and in Theorem~\ref{existence}, we formally
demonstrate that with probability $1$ there exists at least one solution
of the above mentioned linear system of equations, or equivalently, there
exist $P_{0}$-consistent reconstruction operators. Once $P_{0}$-consistent
reconstruction operators are available, local and global
mass matrices are constructed just like the standard MFD method. Therefore,
all the solid foundation of the standard MFD can be readily applied to
our novel curved MFD method.

In Section~\ref{numeric} we test the curved MFD method on grids with highly curved
faces and also curved boundaries, and for all the considered examples,
consistency and the optimal convergence rate are achieved.

We note that the ability of expressing the standard MFD method in a geometric
language is very important for improving its basic building blocks. Thus,
we emphasize the importance of seeing the geometry hidden behind the standard
MFD method. We believe that this geometric viewpoint might also open new
perspectives on the general treatment of curved faces for high-order extensions
of the MFD method like the Virtual Element Method, which is still an open
problem \cite{brezzi2023virtual}. We summarize these observations in
Section~\ref{conclusions}.

\section{Curved grids and geometric entities}
\label{sec2}

Let $\Omega $ be a bounded domain of $\mathbb R^{3}$. Here, by the term
domain, we mean the closure of a connected open subset of
$\mathbb R^{3}$. In particular, $\Omega $ contains its boundary
$\partial \Omega $. A $k$-cell is a $k$-dimensional manifold in
$\mathbb R^{3}$ that is homeomorphic to the closed $k$-dimensional ball.
We equip each $k$-cell in $K$ with an inner orientation
\cite{Tonti2013TheMS}. A \emph{cell complex} $K$ on $\Omega $ is a finite
collection of $k$-cells for $k\in \{0,1,2,3\}$ such that the following
conditions hold \cite{Christiansen2008ACO}:
\begin{enumerate}[(i)]
\item Distinct $k$-cells in $K$ have disjoint interiors.
\item The boundary of any $k$-cell in $K$ is a union of $l$-cells in
$K$ with $l < k$.
\item The union of all cells in $K$ is $\Omega $.
\item The intersection of $k$-cell and an $l$-cell in $K$ with
$l \leq k$ (if not empty) is a union of cells in $K$.
\end{enumerate}

Crucially, $k$-cells as defined above can be \emph{curved}. However, the
standard MFD and DGA methods \cite{cst} require $k$-cells to be
\emph{planar} (or flat). A $k$-cell in $K$ is planar if it is contained
in a $k$-dimensional affine plane of $\mathbb R^{3}$. For our purposes,
$2$-cells in $K$ can be either flat or curved; consequently, also
$1$-cells in $K$ can be either flat, like in a standard polyhedral cell,
or curved. For each 3-cell of $K$ we assume the same geometric regularity
assumption \textbf{(M2)} in \cite{Brezzi2006CONVERGENCEOM}. With this assumption
on $K$, we say that $K$ is a \emph{curved grid} on $\Omega $.

We denote the curved grid $K$ as $K = (N, E, F, C)$, where $N$ the set
of $0$-cells (\emph{nodes}), $E$ the set of $1$-cells (\emph{edges}),
$F$ the set of $2$-cells (\emph{faces}) and $C$ the set of 3-cells (\emph{elements})
in $K$.

We denote by $\abs{\cdot}$ either the cardinality of a set or the measure
of a geometric element in $K$. For example, $\abs{F}$ is the number of
faces in $K$, $\abs{f}$ is the area of face $f \in F$ and $\abs{c}$ is
the volume of element $c \in C$.

Let $X$ be any set among $N,E,F$, $C$ so that elements of $X$ can be
$k$-cells with $k \in \{0,1,2,3\}$. Given a $l$-cell $y$ in $K$, we denote
by $X(y)$ the subset defined by
%
\begin{equation}
X(y) \coloneqq \{ x \in X \mid x \subset \partial y\},
\label{eq1_v1}
\end{equation}
if $k < l$, otherwise,
%
\begin{equation}
X(y) \coloneqq \{x \in X \mid y \subset \partial x\}.
\label{eq2_v1}
\end{equation}
For example, $F(c) = \{ f \in F \mid f \subset \partial c\}$ collects the
faces of $c$ and $C(e) = \{c \in C \mid e \subset \partial c\}$ is the
cluster of elements containing edge $e$.

Let us consider a Cartesian system of coordinates with specified origin
$\bm 0$ and denote by
$\bm x=(x_{1},x_{2},x_{3})^{T} \in \mathbb R^{3}$ the coordinates of its
generic point. Given a vector $\bm x=(x_{1},x_{2},x_{3})^{T}$, we define
${\bm x}|_{(i)}\coloneqq x_{i}$ for $i \in \{1,2,3\}$.

A node $\bm n\in N$ is a point of $\mathbb R^{3}$. A face $f \in F$ is
either a curved or a planar $2$-cell. For each face $f \in F$, we define
its \emph{face vector} $\bm f$ by
%
\begin{equation}
\bm f \coloneqq \int _{f} \hat{\bm n}_{f} \, dS,
\label{face_vector}
\end{equation}
where $\hat{\bm n}_{f}$ denotes the unit vector of $\mathbb R^{3}$ orthogonal
to $f$ (at each point) and oriented according to the right-hand rule with
respect to the orientation of $f$.
%
\begin{remark}
\label{face_vector_same}
The integral (\ref{face_vector}) is the same for every 2-cell spanning the
boundary of $f$ as a consequence of the Divergence Theorem.
\end{remark}

\section{Curved MFD method}
\label{math}

In this section we introduce the basic building blocks of the curved MFD
method. We focus on a curved grid $K=(N,E,F,C)$, and we give specific definitions
only for the vector space $\mathcal F$ of face
\emph{degrees of freedoms} (DoFs) since similar definitions apply to the vector
spaces of DoFs $\mathcal E$ and $\mathcal C$ corresponding to edges and
elements, respectively; see e.g., Section 3 in
\cite{Pitassi2021TheRO}. Overall, the definitions of the vector space
$\mathcal F$ and of the discrete divergence $\mathbb{D}$ and curl operator
$\mathbb{C}$ are the same of the standard MFD method, while the definition
of local reconstruction operators is different and we employ the novel
class of $P_{0}$-consistent local reconstruction operators. Nonetheless, the discrete
inner products $[\cdot , \cdot ]^{\mathcal F}$ and the corresponding derived
differential operators are constructed from $P_{0}$-consistent local reconstruction
operators using the same construction process of the standard MFD method.

\subsection{Projection operators}
\label{sec3.1}

 Let $\bm J$ be a sufficiently regular vector-valued function so that
the integrals of its normal components are well defined on each face
$f \in F$. The \emph{projection operator} $P^{\mathcal F}$ maps
$\bm J$ onto its face DoFs array $\mathbf J = P^{\mathcal F}(\bm J)$, where
the entry ${\mathbf J}|_{(f)}$ of $\mathbf J$ corresponding to face
$f$ is
%
\begin{equation}
{\mathbf J}|_{(f)} \coloneqq \int _{f} \bm J \cdot \hat{\bm n}_{f}\, df.
\label{int_formula}
\end{equation}
The set of all arrays $\mathbf J$ defines the
vector space $\mathcal F= \mathbb R^{\abs{F}}$ of face DoFs.

For each cell $c \in C$, let $\mathcal F_{c}$ be $\abs{F(c)}$-dimensional
vector subspace of $\mathcal F$ obtained by selecting entries corresponding
to each face in $F(c)$. Accordingly, let
$P^{\mathcal F}_{c} : W^{\mathcal F}_{c} \to \mathcal F_{c}$ be the
\emph{local projection operator} (on element $c$), where
$W^{\mathcal F}_{c}$ is a suitable finite-dimensional vector space defined in such a way that $P^{\mathcal F}_{c}$ is bijective on $\mathcal F_{c}$
\cite{Pitassi2021TheRO,da2014mimetic}. We denote by
$\mathbb{P}^{\mathcal F}_{c}$ the matrix associated with the restriction
of $P^{\mathcal F}_{c}$ to the vector subspace of constant vector fields
defined on $c$. $\mathbb{P}^{\mathcal F}$ is a matrix of size
$\abs{F(c)} \times 3$ and reads
%
\begin{equation}
\label{local_projection}
\mathbb{P}^{\mathcal F}_{c} =
\begin{pmatrix}
\vdots
\\
\bm f^{T}
\\
\vdots
\end{pmatrix}
,
\end{equation}
where each $\bm f$ is the face vector of a face $f \in F(c)$. Finally, let
$\mathbb{S}^{\mathcal F}_{c}$ be the \emph{restriction matrix} (on element
$c$) of size $\abs{F(c)} \times \abs{F}$ which maps face DoFs in
$\mathcal F$ to their restriction in $\mathcal F_{c}$. Each entry
${\mathbb{S}^{\mathcal F}_{c}}|_{(f,f')}$, corresponding to a face
$f\in F(c)$ and a face $f' \in F$, is $+1$ if $f = f' \in F(c)$ and
$0$ otherwise.

\subsection{Reconstruction operators}
\label{sec3.2}

Let $R^{\mathcal F}_{c} : \mathcal F_{c} \to W^{\mathcal F}_{c}$ be the
\emph{local reconstruction operator} (on element $c$) defined as the inverse
of map $P^{\mathcal F}_{c}$, which is well-defined thanks to the definition
of the finite-dimensional vector subspace $W^{\mathcal F}_{c}$.
We define the \emph{average} of local reconstruction operator
$R^{\mathcal F}_{c}$ by
%
\begin{equation}
\overline{R^{\mathcal F}_{c}} \coloneqq \frac{1}{\abs{c}} \int _{c} R^{
\mathcal F}_{c} \, dV,
\label{average}
\end{equation}
and we denote by $\mathbb{R}^{\mathcal F}_{c}$ the associated matrix. In
what follows, we will only deal with averages of local reconstruction operators
and for convenience we will drop the words ``average'' and ``local'' and refer
to them simply as reconstruction operators.

Following \cite{Pitassi2021TheRO}, we require the family
$\{\mathbb{R}^{\mathcal F}_{c}\}_{c \in C}$ of reconstruction operators
to satisfy the $P_{0}$-consistency property defined as follows (see also
Definition 4.1 in \cite{Pitassi2021TheRO}).
%
\begin{definition}[$P_{0}$-consistent face reconstruction operators]
\label{consistency}
A collection of face reconstruction operators
$\{\mathbb{R}^{\mathcal F}_{c}\}_{c \in C}$ is said to be $P_{0}$-consistent
if the following two conditions hold:
%
\begin{align}
&\textbf{(P1)} &&\mathbb{R}^{\mathcal F}_{c} \mathbb{P}^{\mathcal F}_{c}
= \mathbb{I}_{3}, \, \forall c \in C,
\label{eq7}\\
&\textbf{(P2)} &&\sum _{f \in F(e)} {\mathbb{C}}|_{(f,e)} \sum _{c
\in C(f)} \abs{c} \, \mathrm{col}_{(f)} {\mathbb{R}^{\mathcal F}_{c}} =
\bm 0, \, \forall e \in E, e \not \subset \partial \Omega .
\label{eq8}
\end{align}
$\mathbb{I}_{3}$ denotes the identity matrix of order 3, while
$\mathbb{C}$ is the face-edge incidence matrix of $K$ defined in Section~\ref{sec3.4} and
$\mathbb{S}_{c}^{\mathcal F}$ is the restriction matrix to element
$c$.
\end{definition}

\begin{remark}
\label{remark_consistency}
\label{observation}
 Note that $P_{0}$-consistency is well-defined on a curved grid $K$. In
fact, it involves only the local projection matrices
$\{\mathbb{P}^{\mathcal F}_{c}\}_{c\in C}$ whose rows are face vectors
of faces of $K$ so that all geometric information of curved faces of
$K$ is solely encoded into their face vectors.
\end{remark}

\begin{figure}
\includegraphics[scale=1.3]{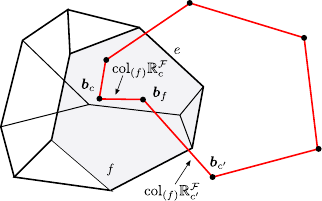}
\caption{Illustration of property \textbf{(P2)} for an internal edge $e$. In
red, the closed path formed by columns (interpreted as vectors of $\mathbb
R^{3}$) of $P_{0}$-consistent reconstruction operators $\{\mathbb{R}^{\mathcal
F}_{c}\}_{c \in C(e)}$. Note how this property links reconstruction operators
defined on different cells (e.g., the two reconstruction operators
$\mathbb{R}^{\mathcal F}_{c}$, $\mathbb{R}^{\mathcal F}_{c'}$ defined on
elements $c, c' \in C(e)$ incident on face $f \in F(e)$)}
\label{P2}
\end{figure}

The definition of $P_{0}$-consistent reconstruction operators boils down
to two requirements. First, property \textbf{(P1)} requires that each local
reconstruction operator $\mathbb{R}^{\mathcal F}_{c}$ is a left inverse
of the local projection operator $\mathbb{P}^{\mathcal F}_{c}$, which is
the accuracy property for element-wise constant vector fields. Second,
it is easy to check that property \textbf{(P2)} has the following geometric
interpretation: for each internal edge $e$ of $K$ (i.e.,
$e \not \subset \partial \Omega $), the vectors
$\abs{c} \, \mathrm{col}_{(f)} {\mathbb{R}^{\mathcal F}_{c}}$ for each
$c \in C(f)$ and each $f \in F(e)$ can be put one after the other to
form a geometric closed path, namely, their vector sum is the zero vector;
see Fig.~\ref{P2} for an illustration.

\subsection{Discrete inner products}
\label{sec3.3}

In order to construct the (face) \emph{mass matrix}
$\mathbb{M}^{\mathcal F}$ associated with the inner product
$[\cdot ,\cdot ]^{\mathcal F}$ on $\mathcal F$, we first construct local
mass matrices $\{\mathbb{M}^{\mathcal F}_{c}\}_{c \in C}$ for each element
$c \in C$ and then we assemble them using a standard FE assembly process.

For each element $c \in C$, let $\mathbb{K}_{c}$ be a symmetric and positive-definite
matrix of order $3$, representing an homogeneous material property in
$c$. We define the local mass matrix $\mathbb{M}^{\mathcal F}_{c}$ (on element $c$) by
%
\begin{equation}
\mathbb{M}^{\mathcal F}_{c} \coloneqq \abs{c}(\mathbb{R}^{\mathcal F}_{c})^{T}
\mathbb{K}_{c} \, \mathbb{R}^{\mathcal F}_{c} + \lambda _{c} \big(
\mathbb{I}_{\abs{F(c)}} - \mathbb{P}^{\mathcal F}_{c}\big((\mathbb{P}^{
\mathcal F}_{c})^{T} \mathbb{P}^{\mathcal F}_{c})^{-1} (\mathbb{P}^{
\mathcal F}_{c})^{T} \big),
\label{local_mass}
\end{equation}
where $\mathbb{I}_{\abs{F(c)}}$ is the identity matrix of order
$\abs{F(c)}$ and $\lambda _{c}$ is the scalar factor
%
\begin{equation}
\lambda _{c} \coloneqq \frac{1}{2} \mathrm{trace} \big( \abs{c} (
\mathbb{R}^{\mathcal F}_{c})^{T} \mathbb{K}_{c} \mathbb{R}^{
\mathcal F}_{c}).
\label{eq10}
\end{equation}
In (\ref{local_mass}), the second term is the so-called
\emph{stabilization part} used in the MFD method (e.g. see \cite[Lemma 2.9 and Corollary 2.1]{Lipnikov2014MimeticFD}),
which ensures that matrix $\mathbb{M}^{\mathcal F}_{c}$ is symmetric, positive-definite
and remains consistent. We note that other expressions for the stabilization
part are possible \cite{Bonelle2015LoworderRO}.

The global mass matrix $\mathbb{M}^{\mathcal F}$ is constructed by assembling
local matrices $\{\mathbb{M}^{\mathcal F}_{c}\}_{c \in C}$ and reads
%
\begin{equation}
\mathbb{M}^{\mathcal F}\coloneqq \sum _{c \in C} (\mathbb{S}^{
\mathcal F}_{c})^{T} \mathbb{M}^{\mathcal F}_{c} \mathbb{S}^{
\mathcal F}_{c}.
\label{global_mass}
\end{equation}

\subsection{Discrete differential operators}
\label{sec3.4}

 Discrete differential operators are represented by incidence matrices
of the grid $K$ just like the standard MFD method.   Therefore, the
\emph{discrete curl operator} is the matrix $\mathbb{C}$ of size
$\abs{F} \times \abs{E}$ of incidence numbers between faces and edges.
Each entry ${\mathbb{C}}|_{(f,e)}$, corresponding to a face $f\in F$ and
an edge $e\in E$, is: $0$, if $e$ is not incident on $f$ (i.e.,
$e \not \subset f$); $+1$, if the orientation of $e$ and $f$ are
\emph{compatible} \cite{Tonti2013TheMS} according to the right-hand rule,
and $-1$ otherwise.   Similarly, the
\emph{discrete divergence operator} is the matrix $\mathbb{D}$ of size
$\abs{C} \times \abs{F}$ of incidence numbers between elements and faces.
Each entry ${\mathbb{D}}|_{(c,f)}$, corresponding to an element
$c \in C$ and a face $f \in F$, is: $0$, if $f$ is not incident on
$c$ (i.e., $f \not \subset c$); $+1$, if the orientation of $f$ and
$c$ are compatible according to the right-hand
rule, and $-1$ otherwise.

 Finally, we can define the
\emph{derived discrete gradient operator} $\tilde{\mathbb{G}}$ as the adjoint
of the discrete divergence operator $\mathbb{D}$ with respect to the inner
product $[\cdot , \cdot ]^{\mathcal F}$ on $\mathcal F$ and  the inner product
$[\cdot , \cdot ]^{\mathcal C}$ on $\mathcal C$ as
%
\begin{equation}
[\tilde{\mathbb{G}}\mathbf U, \,\mathbf J]^{\mathcal F}= [\mathbf U,
\,\mathbb{D} \mathbf J]^{\mathcal C},\qquad \forall \mathbf U \in
\mathcal C,\, \forall \mathbf J \in \mathcal F.
\label{eq12}
\end{equation}
The inner product $[ {\cdot } \, , {\cdot } ]^{\mathcal C}$ on
$\mathcal C$ is defined as in
\cite[Equation (3.12)]{Brezzi2006CONVERGENCEOM} (or equivalently, in
\cite[Equation (2.44)]{Lipnikov2014MimeticFD}); note that its definition
is actually unique since there is only one quadrature rule for the volume
integral on each element that uses a single quadrature point.

\section{Properties of $P_{0}$-consistent reconstruction operators}
\label{curved}

In this section, we derive properties of $P_{0}$-consistent reconstruction
operators. In particular, in Section~\ref{structure_section} we focus on properties
that are valid on general curved grids. Here, the central result is
Theorem~\ref{linear_system_to_local_operators}, which characterizes $P_{0}$-consistent
reconstruction operators as the solution set of the linear system
(\ref{main_linear}). Instead, in Section~\ref{existence_section}, we focus on the practically
important case of curved grids obtained by randomly perturbing internal
nodes of cubic grids, and in {Theorem~\ref{existence} we establish a special existence
result of solutions of linear system {(\ref{main_linear})}, and hence, of
$P_{0}$-consistent reconstruction operators.

\subsection{Structure properties on general curved grids}
\label{structure_section}

To start with, we show that $P_{0}$-consistent reconstruction operators
$\{\mathbb{R}^{\mathcal F}_{c}\}_{c \in C}$ are defined up to an arbitrary
point in $\mathbb R^{3}$; this point can be used to optimize them and we
will exploit this possibility in Section~\ref{numeric}.

\begin{lemma}
\label{node_lemma}
For each $c\in C$, let
$\bm D_{c} \coloneqq (\dots , {\mathbb{D}}|_{(c,f)},\dots )^{T}$ be the
$\abs{F(c)}$-dimensional vector collecting entries
${\mathbb{D}}|_{(c,f)}$ for each $f \in F(c)$. Let
$\{\mathbb{R}^{\mathcal F}_{c}\}_{c \in C}$ be a family of $P_{0}$-consistent
reconstruction operators. Given arbitrary points
$\{\bm b_{c}\}_{c \in C}$ in $\mathbb R^{3}$, we define the family
$\{\mathbb{R}^{\mathcal F}_{\bm b_{c}}\}_{c \in C}$ by
%
\begin{equation}
\mathbb{R}^{\mathcal F}_{\bm b_{c}} \coloneqq \mathbb{R}^{\mathcal F}_{c}
- \frac{1}{\abs{c}} \bm b_{c} \bm D_{c}^{T}, \, \forall c \in C.
\label{virtual_dual_node}
\end{equation}
Then, $\{ \mathbb{R}^{\mathcal F}_{\bm b_{c}} \}_{c \in C}$ is also a family
of $P_{0}$-consistent reconstruction operators.
\begin{proof}
We need to prove that the family
$\{ \mathbb{R}^{\mathcal F}_{\bm b_{c}} \}_{c \in C}$ satisfies properties
\textbf{(P1)} and \textbf{(P2)}.

Proof of \textbf{(P1)}. We note that for each element $c \in C$
%
\begin{equation}
(\mathbb{P}^{\mathcal F}_{c})^{T}\bm D_{c} = \sum _{f \in F(c)} {
\mathbb{D}}|_{(c,f)} \bm f = \sum _{f \in F(c)} \int _{f} {\mathbb{D}}|_{(c,f)}
\bm n_{f} \, dS = \int _{\partial c} \bm n_{\partial c} \, dS = \bm 0,
\label{eq14}
\end{equation}
where the first equality follows from the definition of vector
$\bm D_{c}$ and matrix $\mathbb{P}^{\mathcal F}_{c}$ while the last one follows
from the Divergence Theorem applied to the closed surface
$\partial c$. Therefore, we have that
$\mathbb{R}^{\mathcal F}_{\bm b_{c}} \mathbb{P}^{\mathcal F}_{c} =
\mathbb{R}^{\mathcal F}_{c} \mathbb{P}^{\mathcal F}_{c} = \mathbb{I}_{3}$.

Proof of \textbf{(P2)}. We have that
%
\begin{align}
\begin{split}
\sum _{f \in F(e)} {\mathbb{C}}|_{(f,e)} \sum _{c \in C(f)} \abs{c}
\, \mathrm{col}_{(f)} {\mathbb{R}^{\mathcal F}_{\bm b_{c}}} &=\sum _{f
\in F(e)} {\mathbb{C}}|_{(f,e)} \sum _{c \in C(f)} \abs{c} \,
\mathrm{col}_{(f)} {\mathbb{R}^{\mathcal F}_{c}} - \sum _{f \in F(e)} {
\mathbb{C}}|_{(f,e)} \sum _{c \in C(f)} {\mathbb{D}}|_{(c,f)} \bm b_{c}
\\
&=- \sum _{f \in F(e)} {\mathbb{C}}|_{(f,e)} \sum _{c \in C(f)} {
\mathbb{D}}|_{(c,f)} \bm b_{c},
\end{split}
\label{eq}
\end{align}
because the family $\{\mathbb{R}^{\mathcal F}_{c}\}_{c \in C}$ satisfies
property \textbf{(P2)} by hypothesis. Therefore, we just need to show that
the last term in {(\ref{eq})} vanishes. We have
%
\begin{align}
\begin{split}
\sum _{f \in F(e)} {\mathbb{C}}|_{(f,e)} \sum _{c \in C(f)} {
\mathbb{D}}|_{(c,f)} \bm b_{c} &=\sum _{c \in C(e)} \sum _{f \in F(c)
\cap F(e)} {\mathbb{C}}|_{(f,e)} {\mathbb{D}}|_{(c,f)} \bm b_{c}
\\
&=\sum _{c \in C(e)} ({\mathbb{C}}|_{(f'_{e,c},e)} \mathbb{D}_{f'_{e,c},c}
+ {\mathbb{C}}|_{(f''_{e,c},e)} \mathbb{D}_{f''_{e,c},c}) \bm b_{c} =
\bm 0,
\end{split}
\label{eqq2}
\end{align}
where $F(c) \cap F(e) = \{f'_{e,c}, f''_{e,c}\}$, and the last equality
follows from the identity $\mathbb{D} \mathbb{C} = \bm 0$, due to the cell
complex structure of $K$ (see, e.g.,
\cite[Lemma 3.6]{christiansen2009foundations}).
\end{proof}
\end{lemma}

Recall that the boundary $\partial \Omega $ of $\Omega $ is
\emph{locally flat} if, locally at each of its points, it can be topologically
transformed into a plane of $\mathbb R^{3}$; more precisely, if the following
is true: for each $\bm x\in \partial \Omega $, there exist an open neighborhood
$U_{\bm x}$ of $\bm x$ in $\mathbb R^{3}$ and a homeomorphism
$\phi _{\bm x}:U_{\bm x}\to \mathbb R^{3}$ such that
$\phi _{\bm x}(U_{\bm x}\cap \partial \Omega )=\{(x_{1},x_{2},x_{3})
\in \mathbb R^{3}:x_{3}=0\}$. Observe that, if $\partial \Omega $ is polyhedral
(in the sense of \cite[Remark 9.5, p. 93]{maggi}), then it is also locally
flat. By \cite{brown}, if $\partial \Omega $ is locally flat then
$\partial \Omega $ admits a collar in $\Omega $, i.e., there exists an
open neighborhood of $\partial \Omega $ in $\Omega $, which is homeomorphic
to $\partial \Omega \times [0,1)$. As an immediate consequence, we have:
%
\begin{lemma}
\label{technical_lemma}
If $\partial \Omega $ is locally flat, then $\Omega $ is simply connected
if and only if the interior of $\Omega $ in $\mathbb R^{3}$ is simply connected.
\end{lemma}
For other relevant related results, we refer the reader to \cite{bgf},
especially Corollary 3.5 of that paper.

\begin{lemma}[Structure of $P_{0}$-consistent reconstruction operators]
\label{structure}
Suppose that $\Omega $ is simply connected and has a locally flat boundary.
Let $K=(V,E,F,C)$ be a curved grid on $\Omega $, and let
$\{\mathbb{R}^{\mathcal F}_{c}\}_{c \in C}$ be a family of $P_{0}$-consistent
reconstruction operators. Then, there exists a set of points
$\{\bm b_{c}\}_{c \in C}$ associated with elements of $K$ and a set of
points $\{\bm b_{f}\}_{f \in F}$ associated with faces of $K$ such that
%
\begin{equation}
\abs{c}\, \mathrm{col}_{(f)} {\mathbb{R}^{\mathcal F}_{c}} = {
\mathbb{D}}|_{(c,f)} (\bm b_{f} - \bm b_{c}), \qquad \forall c \in C,
\forall f \in F(c).
\label{eq17}
\end{equation}
\begin{proof}
Let $\tilde{G}=(\tilde{V}, \tilde{E})$ be the graph defined as follows:
the set of vertices $\tilde{V}$ is $\tilde{V} \coloneqq C$; the set of
edges $\tilde{E}$ contains edge $\{c, c'\}$ if and only if
$c \cap c' = f$ for some $f \in F$. We note that graph $\tilde{G}$ is connected
since $K$ is.

Let $\sigma $ be a path in $\tilde{G}$ and denote by
$\tilde{E}_{\sigma }\subset \tilde{E}$ the set of its edges. We fix an
orientation on $\sigma $ by choosing one of the two possible directions
to walk on it. Note that all vertices in a path $\sigma $ are distinct
so that each edge $\{c, c'\} \in \tilde{E}_{\sigma}$ corresponds to a unique
face $f\in F$ such that $c \cap c' = f$. Let
$\bm \sigma \in \mathbb R^{\abs{F}}$ be an array associated to the path
$\sigma $ and defined as follows: ${\bm \sigma }|_{(f)} \neq 0 $ if and
only if $f \in \tilde{E}_{\sigma}$ and, in particular,
${\bm \sigma }|_{(f)} = +1$ if the orientations of $f$ and $\sigma $ (as
fixed above) are consistent with the right-hand rule, and
${\bm \sigma }|_{(f)} = -1$ otherwise. We now show that, if
$\sigma $ is closed path, then property \textbf{(P2)} implies that
%
\begin{equation}
\sum _{f \in \tilde{E}_{\sigma}} {\bm \sigma }|_{(f)} \sum _{c \in C(f)}
\abs{c} \mathrm{col}_{(f)} {\mathbb{R}^{\mathcal F}_{c}} = \sum _{f
\in \tilde{E}_{\sigma}} {\bm \sigma }|_{(f)} \sum _{c \in C(f)} \bm r_{f,c}
= \bm 0,
\label{closed}
\end{equation}
where we have defined
$\bm r_{f,c} \coloneqq \abs{c} \, \mathrm{col}_{(f)} {\mathbb{R}^{
\mathcal F}_{c}}$ just to simplify notation for the rest of the proof.
To prove this fact, let us consider array
$\bm \chi \in \mathbb R^{\abs{E}}$ such that
$\mathbb{C} \bm \chi = \bm \sigma $, which exists thanks to the assumption
that $\Omega $ is simply connected and has a locally flat boundary, so
{Lemma~\ref{technical_lemma}} is verified. It follows that
%
\begin{equation}
\sum _{e \in E, \, {\bm \chi }|_{(e)} \neq 0} {\bm \chi }|_{(e)}
\sum _{f \in F(e)} {\mathbb{C}}|_{(f,e)} \sum _{c \in C(f)} \bm r_{f,c}
=\sum _{f \in \tilde{E}_{\sigma}} {\bm \sigma }|_{(f)} \sum _{c \in C(f)}
\bm r_{f,c} =\bm 0,
\label{eq19}
\end{equation}
where we have used property \textbf{(P2)} and the definition of
$\bm \chi $ that encodes the fact that ``internal'' edges in the support
of array $\bm \chi $ contribute with pairs of opposite sign that mutually
eliminates.

We now fix a vertex $c^{*}$ of $\tilde{G}$. For any other vertex $c$ of
$\tilde{V} \setminus \{c^{*}\}$ we consider a path $\sigma _{c}$ that connects
$c^{*}$ to $c$ in $\tilde{G}$ and we define point
$\bm b_{c,\sigma _{c}}$ (with respect to the fixed origin $\bm 0$) by
%
\begin{equation}
\bm b_{c,\sigma _{c}} \coloneqq \sum _{f' \in \tilde{E}_{\sigma _{c}}}
{\bm \sigma _{c}}|_{(f')} \sum _{c' \in C(f')} \bm r_{f',c'},
\label{formula}
\end{equation}
where we orient the path by walking on it from $c^{*}$ to $c$. Of course,
if there exists a unique path $\sigma _{c}$ connecting $c^{*}$ to
$c$ we set $\bm b_{c} \coloneqq \bm b_{c, \sigma _{c}}$. Otherwise, the
crucial aspect of definition {(\ref{formula})} is that, thanks to
{(\ref{closed})}, the point $\bm b_{c,\sigma _{c}}$ does not depend on the
specific choice of the path $\sigma _{c}$ in $\tilde{G}$ as long as it
connects $c^{*}$ to $c$. As a matter of fact, if we consider a different
path $\sigma '_{c}$ still connecting $c^{*}$ to $c$, then the ``union'' of
these two paths forms the closed path ``$\sigma _{c} \cup \sigma '_{c}$''
for which {(\ref{closed})} and {(\ref{formula})} yields
%
\begin{equation}
\bm b_{c, \sigma _{c}} =\sum _{f' \in \tilde{E}_{\sigma _{c}}} {
\bm \sigma _{c}}|_{(f')} \sum _{c' \in C(f')} \bm r_{f',c'} =\sum _{f'
\in \tilde{E}_{\sigma '_{c}}} {\bm \sigma '_{c}}|_{(f')} \sum _{c'
\in C(f')} \bm r_{f',c'} =\bm b_{c, \sigma '_{c}},
\label{eq21_v1}
\end{equation}
where signs are taken into account to represent the fact that
$\sigma '_{c}$ is walked in opposite sense. Therefore, we set
$\bm b_{c} \coloneqq \bm b_{c,\sigma _{c}} = \bm b_{c,\sigma '_{c}}$.

Next, we define point $\bm b_{f,c}$ by
%
\begin{equation}
\bm b_{f,c} \coloneqq ({\mathbb{D}}|_{(c,f)})^{-1} \bm r_{f,c} +
\bm b_{c}={\mathbb{D}}|_{(c,f)}\bm r_{f,c} + \bm b_{c}
\label{formula2}
\end{equation}
for each $c\in C$ and each $f\in F(c)$. Of course, if $f$ is a boundary
face, then $c$ is the unique element incident on it and we set
$\bm b_{f} \coloneqq \bm b_{f,c}$. Otherwise, $f$ is the common face of
two elements $c,c''$, and we now demonstrate that point
$\bm b_{f,c}$ does not depend on the specific element incident on face
$f$. As a matter of fact, the ``union'' of paths $\sigma _{c}$,
$\sigma _{c''}$, $f$, with $f$ thought as an element of $\tilde{E}$, forms
the closed path ``$\sigma _{c} \cup \sigma _{c''} \cup \{f\}$'' for which
{(\ref{closed})} and {(\ref{formula})}, {(\ref{formula2})} yields:
%
\begin{equation}
\bm b_{c}+{\mathbb{D}}|_{(c,f)}(\bm r_{f,c}+ \bm r_{f,c''})-\bm b_{c''}=
\bm 0
\label{eq23}
\end{equation}
and hence
%
\begin{equation}
\bm b_{f,c} =\bm b_{c}+ {\mathbb{D}}|_{(c,f)} \bm r_{f,c} =\bm b_{c''}-{
\mathbb{D}}|_{(c,f)} \bm r_{f,c''} =\bm b_{c''}+ {\mathbb{D}}|_{(c'',f)}
\bm r_{f,c''} =\bm b_{f,c''}.
\label{eq24}
\end{equation}
Therefore, we set $\bm b_{f} \coloneqq \bm b_{f,c}=\bm b_{f,c''}$.

To conclude, using {(\ref{formula2})} and the definition of vectors
$\bm r_{f,c}$ as given above we have that the collection of all points
$\{\bm b_{f}\}_{f \in F}$ and $\{\bm b_{c}\}_{c \in C}$ satisfy
%
\begin{equation}
\abs{c} \, \mathrm{col}_{(f)} {\mathbb{R}^{\mathcal F}_{c}} = {
\mathbb{D}}|_{(c,f)} (\bm b_{f} - \bm b_{c}),
\label{representation}
\end{equation}
as desired.
\end{proof}
\end{lemma}

The two requirements \textbf{(P1)} and \textbf{(P2)} in
{Definition~\ref{consistency}} can be encoded into a linear system of equations by defining suitable matrices $\mathbb{P}^{\mathcal F}$,
$\mathbb{R}^{\mathcal F}$, $\mathbb{I}$ (depending
on $K$) as follows.
\begin{itemize}
\item Matrix $\mathbb{P}^{\mathcal F}$ of size
$3 \abs{C} \times \abs{F}$. $\mathbb{P}^{\mathcal F}$ is the block matrix
%
\begin{equation}
\mathbb{P}^{\mathcal F}=
\begin{pmatrix}
\vdots
\\
(\mathbb{P}^{\mathcal F}_{c})^{T} (\mathrm{diag} \bm D_{c})\,
\mathbb{S}^{\mathcal F}_{c}
\\
\vdots
\end{pmatrix}
,
\label{def_P}
\end{equation}
where $\mathbb{P}^{\mathcal F}_{c}$ is the local projection matrix, $\mathbb{S}^{\mathcal F}_{c}$ the restriction
matrix and $\mathrm{diag} \bm D_{c}$ is the diagonal matrix of size
$\abs{F(c)} \times \abs{F(c)}$ whose main diagonal is the vector
$\bm D_{c}$ defined in the statement of {Lemma~\ref{node_lemma}}; see
{Fig.~\ref{matrix}} for a graphical representation of matrix
$\mathbb{P}^{\mathcal F}$.

\begin{figure}
\includegraphics[scale=0.8]{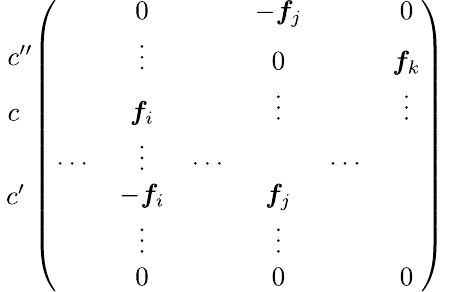}
\caption{Structure of matrix $\mathbb{P}^{\mathcal F}$: face vector of each
internal face appears in two different blocks with opposite signs, where the two
blocks correspond to the unique two elements sharing the internal face; instead,
face vector of each boundary face appears only in one block corresponding to the
unique element containing it. For example, the internal face $f_{i}$ is shared
between elements $c,c'$ and its face vector $\bm f_{i}$ appears with opposite
signs on the rows corresponding to elements $c,c'$; the boundary face
$f_{k}$ is contained in the unique element $c''$ and its face vector $\bm f_{k}$
appears only on the rows corresponding to element $c''$.}
\label{matrix}
\end{figure}
%
\item Matrix $\mathbb{R}^{\mathcal F}$ of size $\abs{F} \times 3$.
$\mathbb{R}^{\mathcal F}$ is the block matrix
%
\begin{equation}
\mathbb{R}^{\mathcal F}=
\begin{pmatrix}
\vdots
\\
\bm b_{f}^{T}
\\
\vdots
\end{pmatrix}
,
\label{generalized_bary}
\end{equation}
where each point $\bm b_{f} \in \mathbb R^{3}$ corresponds to a face
$f \in F$.
\item Matrix $\mathbb{I}$ of size $3\abs{C} \times 3$. $\mathbb{I}$ is
the block matrix
%
\begin{equation}
\mathbb{I} =
\begin{pmatrix}
\vdots
\\
\abs{c} \mathbb{I}_{3}
\\
\vdots
\end{pmatrix}
,
\label{def_I}
\end{equation}
where $\abs{c}$ denotes the volume of an element $c\in C$ and
$\mathbb{I}_{3}$ is the identity matrix of order 3.
\end{itemize}
In what follows, we will also denote matrices
$\mathbb{P}^{\mathcal F}$, $\mathbb{R}^{\mathcal F}$, $\mathbb{I}$ as
$\mathbb{P}^{\mathcal F}_{K}$, $\mathbb{R}^{\mathcal F}_{K}$,
$\mathbb{I}_{K}$ in order to highlight their dependence on geometric elements
of the curved grid $K$.

\begin{remark}[Structure of rows of $\mathbb{P}^{\mathcal F}$]
\label{P_rows}
By definition of block matrix $\mathbb{P}^{\mathcal F}$ in
{(\ref{def_P})}, a row of $\mathbb{P}^{\mathcal F}$ is a row of some block
$(\mathbb{P}^{\mathcal F}_{c})^{T} (\mathrm{diag} \bm D_{c})\,
\mathbb{S}^{\mathcal F}_{c}$ for some element $c\in C$ and thus it collects
one among the three coordinates ${\bm f}|_{(i)}$ with
$i \in \{1,2,3\}$ of face vectors decomposing the boundary of $c$. Therefore,
each row of $\mathbb{P}^{\mathcal F}$ can be put in one-to-one correspondence
with a pair made of an element of $C$ and an index $i \in \{1,2,3\}$, and we can identify the set of rows of $\mathbb{P}^{\mathcal F}$ corresponding
to an index $i \in \{1,2,3\}$ with the set of pairs
$\Pi _{i} \coloneqq C \times \{i\}$. Accordingly, the set of all rows of
$\mathbb{P}^{\mathcal F}$ is the set
$\Pi \coloneqq \bigcup _{i \in \{1,2,3\}} \Pi _{i}$.
\end{remark}

\begin{theorem}[Linear system associated with $P_{0}$-consistent reconstruction operators]
\label{linear_system_to_local_operators}
Let $K=(V,E,F,C)$ be a curved grid on $\Omega $. Let us consider the linear
system of equations
%
\begin{equation}
\mathbb{P}^{\mathcal F}_{K} \mathbb{R}^{\mathcal F}_{K} = \mathbb{I}_{K},
\label{main_linear}
\end{equation}
where $\mathbb{P}^{\mathcal F}_{K}$, $\mathbb{R}^{\mathcal F}_{K}$ and
$\mathbb{I}_{K}$ are defined in {(\ref{def_P})}, {(\ref{generalized_bary})} and
{(\ref{def_I})}, respectively. The following two assertions hold.
\begin{enumerate}[(i)]
\item If $\mathbb{R}^{\mathcal F}_{K}$ is a solution of
{(\ref{main_linear})}, then there exists a family
$\{\mathbb{R}^{\mathcal F}_{c}\}_{c \in C}$ of $P_{0}$-consistent reconstruction
operators.
\item Suppose that $\Omega $ is simply connected and has a locally flat
boundary, and let $\{\mathbb{R}^{\mathcal F}_{c}\}_{c \in C}$ be a family
of $P_{0}$-consistent reconstruction operators. Then, there exists a solution
$\mathbb{R}^{\mathcal F}_{K}$ of {(\ref{main_linear})}.
\end{enumerate}
\begin{proof}
(i) Let $\mathbb{R}^{\mathcal F}$ be a solution of {(\ref{main_linear})}.
We define the family
$\{\mathbb{R}^{\mathcal F}_{\bm 0_{c}}\}_{c \in C }$ by
%
\begin{equation}
\label{formula_op}
(\mathbb{R}^{\mathcal F}_{\bm 0_{c}})^{T} \coloneqq \frac{1}{\abs{c}} (
\mathrm{diag} \bm D_{c})\, \mathbb{S}^{\mathcal F}_{c} \mathbb{R}^{
\mathcal F}, \, \forall c \in C.
\end{equation}
We now show that
$\{\mathbb{R}^{\mathcal F}_{\bm 0_{c}}\}_{c \in C }$ is a family of
$P_{0}$-consistent reconstruction operators, namely, it satisfies properties
\textbf{(P1)} and \textbf{(P2)}.

Proof of \textbf{(P1)}. By definition of block matrix
$\mathbb{P}^{\mathcal F}$ in {(\ref{def_P})}, we have that equation
$(\mathbb{P}^{\mathcal F}_{c})^{T} (\mathrm{diag} \bm D_{c})\,
\mathbb{S}^{\mathcal F}_{c} \mathbb{R}^{\mathcal F}= \abs{c}
\mathbb{I}_{3}$ holds for each element $c \in C$. Transposing members and
using the definition of $\mathbb{R}^{\mathcal F}_{\bm 0_{c}}$ we get
$\mathbb{R}^{\mathcal F}_{\bm 0_{c}} \mathbb{P}^{\mathcal F}_{c} =
\mathbb{I}_{3}$.

Proof of \textbf{(P2)}. We have that
%
\begin{equation}
\sum _{f \in F(e)} {\mathbb{C}}|_{(f,e)} \sum _{c \in C(f)} \abs{c}
\, \mathrm{col}_{(f)} {\mathbb{R}^{\mathcal F}_{\bm 0_{c}}} =\sum _{f
\in F(e)} {\mathbb{C}}|_{(f,e)} \sum _{c \in C(f)} {\mathbb{D}}|_{(c,f)}
\bm b_{f}=\bm 0
\label{eq3}
\end{equation}
because the second term in {(\ref{eq3})} vanishes by reasoning as in
{(\ref{eqq2})} in the proof of {Lemma~\ref{node_lemma}}.

(ii) Let $\{\mathbb{R}^{\mathcal F}_{c}\}_{c \in C}$ be a family of
$P_{0}$-consistent reconstruction operators and let us apply
{Lemma~\ref{structure}} to it. Hence, there exist points
$\{\bm b_{c}\}_{c \in C}$ and $\{\bm b_{f}\}_{f \in F}$ such that
$\abs{c}\, \mathrm{col}_{(f)} {\mathbb{R}^{\mathcal F}_{c}} ={
\mathbb{D}}|_{(c,f)}(\bm b_{f} - \bm b_{c})$ for each $c \in C$ and each
$f \in F(c)$. Let
$\{\mathbb{R}^{\mathcal F}_{\bm 0_{c}}\}_{c \in C}$ be the family of
$P_{0}$-consistent reconstruction operators obtained by applying
{Lemma~\ref{node_lemma}} to $\{\mathbb{R}^{\mathcal F}_{c}\}_{c \in C}$\vspace{3pt} with
points $\{\bm b_{c}\}_{c \in C}$ up to a minus sign; it is then clear that
$\abs{c}\, \mathrm{col}_{(f)} {\mathbb{R}^{\mathcal F}_{\bm 0_{c}}} = {
\mathbb{D}}|_{(c,f)} \bm b_{f}$. Now consider equations
$\abs{c}\mathbb{R}^{\mathcal F}_{\bm 0_{c}} \mathbb{P}^{\mathcal F}_{c}
= \abs{c} \mathbb{I}_{3}$ for each element $c \in C$. By transposing members,
using the definition of block matrices $\mathbb{P}^{\mathcal F}$,
$\mathbb{I}$ in {(\ref{def_P})}, {(\ref{def_I})}, and using the points
$\{\bm b_{f}\}_{f \in F}$ to define $\mathbb{R}^{\mathcal F}$ as in
{(\ref{generalized_bary})}, we have that the matrix
$\mathbb{R}^{\mathcal F}$ is a solution of linear system
{(\ref{main_linear})}.
\end{proof}
\end{theorem}

\begin{remark}[Dual grid structure associated with $P_{0}$-consistent reconstruction operators]
\label{geometry_consistency}
Corresponding to $P_{0}$-consistent reconstruction operators solution of
linear system {(\ref{main_linear})} there is a ``generalized'' dual grid structure
as follows. We think each point $\bm b_{f}$ in
{(\ref{generalized_bary})} as a ``generalized barycenter'' of face $f$ (hence,
the notation $\bm b_{f}$), and we think each point $\bm b_{c}$ as a ``generalized
barycenter'' of element $c$. Next, we think each pair
$(\bm b_{f}, \bm b_{c})$, where $f$ is a face of element $c$, as a ``generalized
dual edge'' with suitable orientation given by the pair order and depending
on the coefficient
${\mathbb{D}}|_{(c,f)}$. See Section~\ref{numeric} for instances of such dual
grid structures constructed for the considered test cases.
\end{remark}

Let $T$ be a rigid translation of $\mathbb R^{3}$ and denote by
$T(K)$ the curved grid obtained by applying $T$ to each node, edge, face
and cell of $K$. Then, it is easy to see that corresponding linear system
$\mathbb{P}^{\mathcal F}_{T(K)} \mathbb{R}^{\mathcal F}_{T(K)} =
\mathbb{I}_{T(K)}$ is equal to {(\ref{main_linear})}, as meaning that
$\mathbb{P}^{\mathcal F}_{T(K)} = \mathbb{P}^{\mathcal F}_{K}$ and
$\mathbb{I}_{T(K)}=\mathbb{I}_{K}$. The next result discusses the case
of rigid rotations of $K$.
%
\begin{proposition}[Transformation of linear system {(\ref{main_linear})} under rigid rotations of $K$]
\label{rotation_lemma}
Let $L$ be a rigid rotation of $\mathbb R^{3}$ about the origin
$\bm 0$ and let $\mathbb{L}$ the corresponding matrix. Denote by
$L(K)$ the curved grid obtained by applying $L$ to each node, edge, face
and cell of $K$. Denote by $\mathbb{P}^{\mathcal F}_{L(K)}$,
$\mathbb{R}^{\mathcal F}_{L(K)}$ and $\mathbb{I}_{L(K)}$ the matrices defined
via elements of $L(K)$ as in {(\ref{def_P})}, {(\ref{generalized_bary})} and
{(\ref{def_I})}, respectively. Accordingly, let us consider linear system
%
\begin{equation}
\mathbb{P}^{\mathcal F}_{L(K)} \mathbb{R}^{\mathcal F}_{L(K)} =
\mathbb{I}_{L(K)},
\label{rotated_main_linear}
\end{equation}
and define
%
\begin{align}
\mathbb{R}^{\mathcal F}_{K} := \mathbb{R}^{\mathcal F}_{L(K)}
\mathbb{L}.
\label{law}
\end{align}
Then, $\mathbb{R}^{\mathcal F}_{K}$ in {(\ref{law})} is a solution of linear
system {(\ref{main_linear})} if and only if
$\mathbb{R}^{\mathcal F}_{L(K)}$ is a solution of linear system
{(\ref{rotated_main_linear})}.
\begin{proof}
Let $f$ be face of $K$ and denote by $f'=L(f)$ be the corresponding rotated
face in $L(K)$. Hence, we have that $\bm f' = \mathbb{L} \bm f$. Using the block structure
of matrix $\mathbb{P}^{\mathcal F}_{K}$ as in {(\ref{def_P})}, we infer that
$\mathbb{P}^{\mathcal F}_{L(K)} = (\mathbb{I}_{\,\abs{C}} \otimes
\mathbb{L}) \mathbb{P}^{\mathcal F}_{K}$, where
$\mathbb{I}_{\,\abs{C}}$ is the identity matrix of order $\abs{C}$. Next, using the block structure
of matrix $\mathbb{I}_{K}$ as in {(\ref{def_I})},
we have that $\mathbb{I}_{L(K)} = \mathbb{I}_{K}$ since $L$ leave volumes
of rotated elements invariant.

Consider now a solution $\mathbb{R}^{\mathcal F}_{L(K)}$ of linear system
{(\ref{rotated_main_linear})}. Then, by multiplying
{(\ref{rotated_main_linear})} by
$(\mathbb{I}_{\,\abs{C}} \otimes \mathbb{L})$ on the left, and by
$\mathbb{L}$ on the right, we get
%
\begin{align}
\label{reasoning}
\begin{split}
(\mathbb{I}_{\,\abs{C}} \otimes \mathbb{L}^{T})\, ( \mathbb{P}^{
\mathcal F}_{L(K)} \mathbb{R}^{\mathcal F}_{L(K)} ) \,\mathbb{L} &= (
\mathbb{I}_{\,\abs{C}} \otimes \mathbb{L}^{T})\, (\mathbb{I}_{L(K)})
\, \mathbb{L} \Rightarrow
\\
\big((\mathbb{I}_{\,\abs{C}} \otimes \mathbb{L}^{T}) (\mathbb{I}_{\,
\abs{C}} \otimes \mathbb{L}) \big) \,\mathbb{P}^{\mathcal F}_{K} \, (
\mathbb{R}^{\mathcal F}_{L(K)} \mathbb{L}) &= \mathbb{I}_{L(K)}
\Rightarrow
\\
\mathbb{P}^{\mathcal F}_{K} \, ( \mathbb{R}^{\mathcal F}_{L(K)}
\mathbb{L}) &= \mathbb{I}_{K},
\end{split}
\end{align}
where we have used the expressions of rotated matrices
$\mathbb{P}^{\mathcal F}_{L(K)}, \mathbb{I}_{L(K)}$ together with the block
structure of matrix $\mathbb{I}_{\,\abs{C}} \otimes \mathbb{L}$ and the
fact that $\mathbb{L}$ is orthogonal since is the matrix associated with
the rigid rotation $L$ of $\mathbb R^{3}$. The last relation in
{(\ref{reasoning})} shows that $\mathbb{R}^{\mathcal F}_{K}$ in
{(\ref{law})} is a solution of {(\ref{main_linear})}, and the converse result
follows immediately since $L^{-1}$ is again a rotation of
$\mathbb R^{3}$.
\end{proof}
\end{proposition}

\subsection{Existence results on cubic grids with randomly perturbed nodes}
\label{existence_section}

Let $K=(V,E,F,C)$ be a curved grid on $\Omega $. We now want to determine
conditions on the geometry of $K$ for which linear system
{(\ref{main_linear})} is compatible. Of course such conditions should restrict
the class of curved grids since matrices
$\mathbb{P}^{\mathcal F}_{K}, \mathbb{R}^{\mathcal F}_{K}$,
$\mathbb{I}_{K}$ in {(\ref{main_linear})} depend on the specific geometry
of curved grid $K$. Therefore, we wish to solve the following general problem:
\begin{description}
\item[(Pr1)] Characterize the set $\mathfrak C$ of all curved grids
$K$ for which linear system {(\ref{main_linear})} is compatible.
\end{description}
A specific version of \textbf{(Pr1)} is the following:
\begin{description}
\item[(Pr2)] Characterize the subset
$\mathfrak R \subset \mathfrak C$ of all curved grids $K$ such that matrix
$\mathbb{P}^{\mathcal F}_{K}$ has full rank, i.e.,
$\rank \mathbb{P}^{\mathcal F}_{K} = 3 \abs{C}$.
\end{description}

Problem \textbf{(Pr1)} seems to be very difficult in general and we do not
tackle it in the present contribution. Instead, we focus on problem
\textbf{(Pr2)} and we determine a sufficient and a necessary condition for
the curved grid $K$ to be contained in the class $\mathfrak R$.

\subsubsection{A necessary condition}
\label{sec4.2.1}

To derive a necessary condition, it is clear that
$\rank \mathbb{P}^{\mathcal F}_{K} = 3\abs{C}$ implies
$\abs{F} \geq 3 \abs{C}$. We now discuss the implications of the condition
$\abs{F} \geq 3 \abs{C}$ on the geometry of curved grids in
$\mathfrak R$. Let us assume that each element $c\in C$ has the same number
$k$ of faces.
Consider the set $Q$ of all pairs $(f,c)$ for each
$c \in C$ and each $f \in F(c)$. Since every element has $k$ faces, we
have $\abs{Q}=k \abs{C}$. Moreover, since every internal face is shared
by exactly two elements and every boundary face is incident to exactly
one element, we have $\abs{Q} = 2 \abs{F_{int}} + \abs{F_{b}}$, where sets
$F_{int}$ and $F_{b}$ partition $F$ into internal and boundary faces, respectively.
By combining these equations and plugging into
$\abs{F} \geq 3 \abs{C}$, we get the equivalent condition
$(k-6)\abs{F_{int}} \geq (3-k) \abs{F_{b}}$, which holds independently
of $\abs{F}$ for $k \geq 6$. In order to practically deal with the cases
$k=4$ or $k=5$ we propose the following construction. For each element
$c \in C$ we select a face $f \in F(c)$ to form a sequence of pairs
$Q'=(\dots , (f,c), \dots )$ such that each face $f$ does not appear in more
than one pair in $Q'$. Then, if $k=4$, we partition each face appearing
in a pair of the sequence $Q'$ into three polygons, and we deform them
to form three linearly independent face vectors. Similarly, if $k=5$, we
partition each face appearing in a pair of the sequence $Q'$ into two polygons,
and we deform them to form two linearly independent face vectors.
In both
cases, in the end we get a curved grid whose elements satisfy
the condition $k\geq 6$. Note that a sequence $Q'$ as above can be constructed
efficiently in linear time and space using standard spanning tree constructions
(e.g. see \cite{Pitassi2021InvertingTD} for a detailed description of such
a construction). The price to pay is the additional number of face DoFs
which shifts from $\abs{F}$ to $\abs{F}+2 \abs{C}$ for $k=4$ and $\abs{F}+ \abs{C}$ for $k=5$.

\subsubsection{A sufficient condition}
\label{sufficient_condition}

To derive a sufficient condition, we restrict to a practically important
class of curved grids and for it we prove that the condition
$\abs{F} \geq 3 \abs{C}$ is sufficient in a very precise sense. Let
${K_{h}}=(N_{h},E_{h},F_{h},C_{h})$ be a cubic grid whose elements in
$C_{h}$ are cubes of edge length $h$ that are, up to a rigid rotation,
subsets of $\mathbb R^{3}$ of the form $\bm z + [0,h]^{3}$ for some
$\bm z \in \mathbb R^{3}$. Following \cite{brezzi2007new}, we now define
the class $\mathfrak M({K_{h}})$ of curved grids constructed from
${K_{h}}$ as follows. Let $P_{h}$ be the cube $[-0.4h, 0.4h]^{3}$, let
$N_{h,int} \subset N_{h}$ be the set of internal nodes of ${K_{h}}$, and
let $P:=\prod _{n\in N_{h,int}}P_{h}=(P_{h})^{\abs{N_{h,int}}}$. Note that
$P$ is a subset of $\mathbb R^{3\abs{N_{h,int}}}$ whose interior is the
non-empty open set $\prod _{n\in N_{h,int}}(-0.4h,0.4h)^{3}$. For each
$\bm n \in N_{h,int}$, consider a point $\bm \xi _{\bm n}$ in
$P_{h}$, and denote by
$\Xi = (\dots , \bm \xi _{\bm n} , \dots )^{T}$ a $3\abs{N_{h,int}}$-dimensional
array in $P$. For each array $\Xi \in P$, we define the curved grid
$K_{\Xi}=(N_{\Xi}, E_{\Xi}, F_{\Xi}, C_{\Xi})$ associated with
$\Xi $ as follows:
\begin{itemize}
\item Set of nodes $N_{\Xi}$. For each node
$\bm n \in N_{h}$ corresponds the \emph{moved} node
$\bm n_{\Xi }\in N_{h,\Xi}$ defined by
%
\begin{equation}
\bm n_{\Xi} \coloneqq
\begin{cases}
\bm n + \bm \xi _{\bm n} & \text{ if } \bm n \in N_{h,int},\\
\bm n & \text{ if } \bm n \in N_{h,b}.
\end{cases}
\label{movement}
\end{equation}
In other words, each internal node $\bm n \in N_{h,int}$ is shifted in a cubic region
of edge length $0.8h$ that is centered at point $\bm n$ and whose edges
are aligned with the axes of cubic grid ${K_{h}}$. Instead, each boundary
node $\bm n \in N_{h,b}$ is left invariant.
%
\item Set of edges $E_{\Xi}$. For each edge $e \in E_{h}$ with boundary
nodes $\{\bm n_{1}, \bm n_{2}\} = N(e)$ corresponds the edge (segment)
$e_{\Xi}$ joining the moved nodes $\bm n_{1,\Xi}, \bm n_{2,\Xi}$.
\item Set of faces $F_{\Xi}$. For each face $f \in F_{h}$ with boundary
nodes $\{\bm n_{1}, \bm n_{2}, \bm n_{3}, \bm n_{4}\} = N(f)$ and such
that $\partial f$ is equal to the following sum of oriented segments
$\partial f = [\bm n_{1}, \bm n_{2}] + [\bm n_{2}, \bm n_{3}] + [
\bm n_{3}, \bm n_{4}] + [\bm n_{4}, \bm n_{1}]$ according to a fixed (arbitrary)
orientation, corresponds the \emph{piecewise linear} face $f_{\Xi}$ made
of two triangles $t_{f}^{(1)}$ and $t_{f}^{(2)}$ with vertices
$\{\bm n_{1,\Xi}, \bm n_{2,\Xi}, \bm n_{3,\Xi}\}$ and
$\{\bm n_{3,\Xi}, \bm n_{4,\Xi}, \bm n_{1,\Xi}\}$, respectively, and oriented
so that
$\partial t_{f}^{(1)} = [\bm n_{1,\Xi}, \bm n_{2,\Xi}] + [\bm n_{2,
\Xi}, \bm n_{3,\Xi}] + [\bm n_{3,\Xi}, \bm n_{1,\Xi}]$ and
$\partial t_{f}^{(2)} = [\bm n_{1,\Xi}, \bm n_{3,\Xi}] + [\bm n_{3,\Xi}, \bm n_{4,\Xi}] + [\bm n_{4,
\Xi}, \bm n_{1,\Xi}]$. Therefore, thanks
to the movement {(\ref{movement})}, the boundary of each face of
$f_{\Xi}$ is not necessarily a planar quadrilateral; see
{Fig.~\ref{curved_cube}} for an illustration.
\item Set of cells $C_{\Xi}$. For each cell $c \in C_{h}$ corresponds the
cell $c_{\Xi}$ whose boundary is the union of moved faces $f_{\Xi}$ with
$f \in F(c)$.
\end{itemize}
 We define $\mathfrak M({K_{h}})$ as the set of all curved grids
$K_{\Xi}$ for $\Xi \in P$.

\begin{figure}
\includegraphics[scale=1.3]{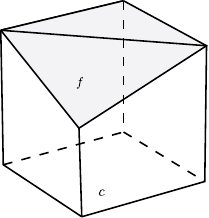}
\caption{A cube $c$ with a curved face $f$ that is a piecewise linear surface
made of two triangles.}
\label{curved_cube}
\end{figure}

The next preliminary result shows that cubic grids ${K_{h}}$ belong to the class
$\mathfrak R$.
%
\begin{lemma}
\label{rank_theorem}
Let ${K_{h}}=(V_{h},E_{h},F_{h},C_{h})$ be a cubic grid. Then, the matrix
$\mathbb{P}^{\mathcal F}_{{K_{h}}}$ has full rank, i.e.,
$\rank \mathbb{P}^{\mathcal F}_{{K_{h}}}= 3 \abs{C_{h}}$. In particular,
$\abs{F_{h}} \geq 3 \abs{C_{h}}$.
\begin{proof}
First, observe that the axes of the cubic grid ${K_{h}}$ are aligned with
the coordinate axes of the chosen Cartesian system of coordinates. For
the rest of proof, we denote $\mathbb{P}^{\mathcal F}_{{K_{h}}}$ by
$\mathbb{P}^{\mathcal F}$, where we drop subscript ${K_{h}}$ just to simplify
notation.

Condition $\rank \mathbb{P}^{\mathcal F}= 3 \abs{C_{h}}$ holds if and only
if there is no linear dependent row vectors of
$\mathbb{P}^{\mathcal F}$. Let us thus assume that, for the sake of contradiction,
that the row vector $\mathrm{row}_{(c,i)} {\mathbb{P}^{\mathcal F}}$, corresponding
to element $c \in C$ and $i$-th coordinate for some
$i \in \{1,2,3\}$, can be written as linear combination of other rows of
$\mathbb{P}^{\mathcal F}$, namely, there exist real coefficients
$\{p_{c',j}\}_{(c',j) \in \Pi \setminus \{(c,i)\}}$ such that
%
\begin{equation}
\mathrm{row}_{(c,i)} {\mathbb{P}^{\mathcal F}} + \sum _{(c',j) \in
\Pi \setminus \{(c,i)\}} p_{c',j}\, \mathrm{row}_{(c',j)} {\mathbb{P}^{
\mathcal F}} = \bm 0.
\label{expression_temp}
\end{equation}
Here $\Pi $ denotes the set of indices defined earlier in
{Remark~\ref{P_rows}}. Let $v$ be the $\abs{F}$-dimensional vector that appears
in the first member of {(\ref{expression_temp})}. Without loss of generality,
we can assume that $i=1$ and $h=1$. Let $f_{0}$ be the face of $c$ such
that $\bm f_{0}=(1,0,0)^{T}$. Note that $f_{0}$ is an internal face of
$K_{h}$. Otherwise, the column vector
$\mathrm{col}_{(f_{0})} {\mathbb{P}^{\mathcal F}}$ has exactly one non-zero
entry ${\mathbb{P}^{\mathcal F}}|_{(c,f_{0})}=1$ corresponding to element
$c$, and the entry ${v}|_{(f_{0})}$ would be $1$ instead of $0$ as imposed
by {(\ref{expression_temp})}. Hence $f_{0}$ is an internal face of
$K_{h}$ and there exists, and is unique, a cube $c_{1}$ of $K_{h}$ such
that $c \cap c_{1} = f_{0}$. Now, the column vector
$\mathrm{col}_{(f_{0})} {\mathbb{P}^{\mathcal F}}$ has exactly two non-zero
entries ${\mathbb{P}^{\mathcal F}}|_{(c,f_{0})}=1$,
${\mathbb{P}^{\mathcal F}}|_{(c_{1},f_{0})}=-1$ corresponding to $c$,
$c_{1}$. According to {(\ref{expression_temp})}, it holds that
$p_{c_{1},1}=1$. By repeating the above reasoning and using the fact that
$K_{h}$ is bounded, we obtain a sequence
$c = c_{0}, c_{1}, \dots , c_{n}$ such that
$c_{l} \cap c_{l+1} = {f_{l}}$ with ${\bm f_{l}}|_{(1)} \neq 0$ for
$l \in \{0, \dots , n-1\}$, $f_{l}$ is an internal face of $K_{h}$ for
$l \in \{0, \dots , n-1\}$ and the face $f_{n}$ of $c_{n}$ opposite to
$f_{n-1}$ is a boundary face of $K_{h}$. It follows that the column vector
$\mathrm{col}_{(f_{n})} {\mathbb{P}^{\mathcal F}}$ has exactly one non-zero
entry ${\mathbb{P}^{\mathcal F}}|_{(c_{n},f_{n})}=1$. As a consequence,
we should have that ${v}|_{(f_{n})}=1$ and not ${v}|_{(f_{n})}=0$ as imposed
by {(\ref{expression_temp})}. This gives the desired contradiction.
\end{proof}
\end{lemma}

In order to state {Theorem~\ref{existence}}, we first recall some basic facts about
\emph{real algebraic subsets} of $\mathbb R^{m}$ with $m$ an arbitrary natural
number \cite{bochnak2013real}. A subset of $\mathbb R^{m}$ is called
\textit{real algebraic} if it is the zero set of a family of polynomials
with (real coefficients and) variables in $\mathbb R^{m}$. A basic result
is that every real algebraic subset $X$ of $\mathbb R^{m}$ can be written
as the union of finite family $\{M_{i}\}_{i}$ of smooth submanifolds
$M_{i}$ of $\mathbb R^{m}$. The dimension $\dim (X)$ of $X$ can be defined
as the maximum dimension of the $M_{i}$'s. If $\dim (X)<m$ then each
$M_{i}$ has dimension $<m$ and hence each $M_{i}$ and the whole $X$ have
measure zero in $\mathbb R^{m}$. In particular, the interiors in
$\mathbb R^{m}$ of each $M_{i}$ and of $X$ are empty. It is important to
observe that, if $X$ is properly contained in $\mathbb R^{m}$, i.e.,
$X\neq \mathbb R^{m}$, then $\dim (X)<m$ so $X$ has measure zero in
$\mathbb R^{m}$. See \cite[Chapters 2, 3 and 9]{bochnak2013real} and
\cite[Chapter 1, Section 0 and Chapter 3, Section 1]{hirsch2012differential}.

\begin{theorem}
\label{existence}
Let ${K_{h}}=(N_{h}, E_{h}, F_{h},C_{h})$ be a cubic grid and let
$K_{\Xi}=(N_{\Xi}, E_{\Xi}, F_{\Xi}, C_{\Xi})$ be the curved grid in
$\mathfrak M(K_{h})$ obtained by the $\Xi $-movement {(\ref{movement})} of
nodes of $K_{h}$. Let $P'$ be the set of all
$\Xi \in P \subset \mathbb R^{3\abs{N_{h,int}}}$ such that the matrix
$\mathbb{P}^{\mathcal F}_{K_{\Xi}}$ has not full rank, namely,
%
\begin{equation}
P':=\big\{\Xi \in P:\mathrm{rank}(\mathbb{P}^{\mathcal F}_{K_{\Xi}})<3
\abs{C_{h}}\}.
\label{singular}
\end{equation}
Then, $P'$ has measure zero in $\mathbb R^{3\abs{N_{h,int}}}$.
\begin{proof}
Let $f_{\Xi}\in F_{\Xi}$. Rearranging indices if necessary, we can assume
that $\partial f_{\Xi}$ is equal to the following sum of oriented segments
$\partial f_{\Xi}=[\bm n_{1,\Xi},\bm n_{2,\Xi}]+[\bm n_{2,\Xi},\bm n_{3,
\Xi}]+[\bm n_{3,\Xi},\bm n_{4,\Xi}]+[\bm n_{4,\Xi},\bm n_{1,\Xi}]$. As
a consequence, we have
%
\begin{align}
\begin{split}
\boldsymbol{f}_{\Xi}=&\frac{1}{2}\big(\bm n_{2,\Xi}-\bm n_{1,\Xi}
\big)\times \big(\bm n_{3,\Xi}- \bm n_{1,\Xi} \big)
\\
&+\frac{1}{2}\big(\bm n_{3,\Xi}-\bm n_{1,\Xi})\times \big(\bm n_{4,
\Xi} - \bm n_{1,\Xi} \big),
\end{split}
\label{face_vector_form}
\end{align}
where ``$\times $'' indicates the vector product in $\mathbb R^{3}$. Using
the expression {(\ref{movement})} of moved nodes of $K_{\Xi}$, we infer that
each of the three components of the vector $\boldsymbol{f}_{\Xi}$ has a
polynomial expression in the variables
$(\xi _{\bm n_{1},1},\xi _{\bm n_{1},2},\xi _{\bm n_{1},3},\ldots ,
\xi _{\bm n_{4},1},\xi _{\bm n_{4},2},\xi _{\bm n_{4},3})^{T} = (
\bm \xi _{\bm n_{1}}, \bm \xi _{\bm n_{2}}, \bm \xi _{\bm n_{3}},
\bm \xi _{\bm n_{4}})^{T} \in (\mathbb R^{3})^{4}$. Note that the latter
assertion remains true even if a node $\bm n_{j,\Xi}=\bm n_{j}$ of
$f_{\Xi}$ is a boundary node of $K_{h}$ and thus is not moved by transformation
{(\ref{movement})}. Indeed, in this case the mentioned polynomial expressions
are constant in the variables
$(\xi _{\bm n_{j},1},\xi _{\bm n_{j},2},\xi _{\bm n_{j},3})^{T}=
\bm \xi _{\bm n_{j}}$.

Recall that $3\abs{C_{h}} \leq |F_{h}|$ by {Lemma~\ref{rank_theorem}}. Denote
by $S$ the family of all subsets of $\{1,\ldots ,|F_{h}|\}$ having cardinality
$3\abs{C_{h}}$. For each $F_{s}\in S$, denote by
${\mathbb{P}^{\mathcal F}_{K_{\Xi }}}|_{(\Pi \times F_{s})}$ the
$3 \abs{C_{h}} \times 3 \abs{C_{h}}$-submatrix of
$\mathbb{P}^{\mathcal F}_{K_{\Xi}}$ obtained selecting all the
$3\abs{C_{h}}$ rows and the $3\abs{C_{h}}$ columns of
$\mathbb{P}^{\mathcal F}_{K_{\Xi}}$ whose indices are the elements of the
set $F_{s}$. Thus,
$\{{\mathbb{P}^{\mathcal F}_{K_{\Xi }}}|_{(\Pi \times F_{s})}\}_{F_{s}
\in S}$ is the family of all $3\abs{C_{h}}\times 3\abs{C_{h}}$-submatrices
of $\mathbb{P}^{\mathcal F}_{K_{\Xi}}$. Define the function
$g: P \to \mathbb R$ by setting
%
\begin{equation}
g(\Xi ):=\sum _{F_{s} \in S} \Big(\det \big({\mathbb{P}^{\mathcal F}_{K_{
\Xi }}}|_{(\Pi \times F_{s})}\big)\Big)^{2},
\label{poly}
\end{equation}
where ``$\det $'' is the determinant operation.

Note that $g$ is a polynomial function in the variables
$(\dots , \xi _{\bm n,1},\xi _{\bm n,2},\xi _{\bm n,3}, \dots )^{T}=(
\dots , \bm \xi _{\bm n}, \dots )^{T} = \Xi \in \mathbb R^{3
\abs{N_{h,int}}}$ since each entry of the matrix
$\mathbb{P}^{\mathcal F}_{K_{\Xi}}$ is. Observe that $g(\Xi )=0$ if and
only if
$\det ({\mathbb{P}^{\mathcal F}_{K_{\Xi }}}|_{(\Pi \times F_{s})})=0$ for
all subsets $F_{s}\in S$, which is in turn equivalent to assert
$\rank \mathbb{P}^{\mathcal F}_{K_{\Xi}} <3\abs{C_{h}}$. This proves that
$P'=P\cap X$, where $X$ is the real algebraic subset of
$\mathbb R^{3 \abs{N_{h,int}}}$ defined by
$X:=\{\Xi \in \mathbb{\mathbb R}^{3 \abs{N_{h,int}}}:g(\Xi )=0\}$.

It is important to observe that the polynomial $g$ in {(\ref{poly})} is non-zero,
i.e., at least one coefficient in the expression of the polynomial
$g$ as linear combination of monomials is different from zero. Indeed,
by {Lemma~\ref{rank_theorem}}, the rank of
$\mathbb{P}^{\mathcal F}_{K_{h}}$ is $3\abs{C_{h}}$. On the other hand,
if $O\in P$ denotes the point associated with the zero vector of
$\mathbb R^{3 \abs{N_{h,int}}}$, then $K_{O}=K_{h}$, and it follows that
$g(O)\neq 0$, i.e., $O\notin X$. Thus,
$X\neq \mathbb R^{3 \abs{N_{h,int}}}$ and hence $X$ has measure zero in
$\mathbb R^{3 \abs{N_{h,int}}}$. Evidently, the same is true for the subset
$P'$ of $X$.
\end{proof}
\end{theorem}

The practical meaning of {Theorem~\ref{existence}} is as follows. Let us consider
random values for
$\Xi = (\dots , \xi _{\bm n,1}, \xi _{\bm n,2},\xi _{\bm n,3},\dots )^{T}
\in P$ where each number $\xi _{\bm n,i}$ for $\bm n \in N_{h,int}$,
$i \in \{1,2,3\}$ is drawn independently from a uniform distribution on
$[-0.4h,0.4h]$. In this case, it may happen that we hit a point
$\Xi \in P'$ so that
$\rank \mathbb{P}^{\mathcal F}_{K_{\Xi}}<3\abs{C_{h}}$. However, thanks
to {Theorem~\ref{existence}}, the probability of this to happen is $0$.

We note that the reasoning behind {Theorem~\ref{existence}} can be extended to other
families of polyhedral grids with randomly perturbed nodes: it is just
sufficient to check whether matrix $\mathbb{P}^{\mathcal F}_{K}$ has full
rank for a given polyhedral grid $K$.

\subsection{Link between $P_{0}$-consistent and admissible reconstruction operators}
\label{comparison}

The standard MFD method defines the class of \emph{admissible} local reconstruction
operators $\{R^{\mathcal F}_{c}\}_{c \in C}$ that satisfy five properties
\textbf{(R1)-(R5)} as defined in \cite{da2014mimetic,Brezzi2014MimeticSP}.
The next result illuminates the link between admissible reconstructions
operators and $P_{0}$-consistent reconstruction operators.
%
\begin{proposition}
\label{fundamental}
Let $K=(N,E,F,C)$ be either a polyhedral or a curved grid. For each element
$c \in C$, define the \emph{dual edge} of face $f \in F(c)$ (restricted
to $c$) by
%
\begin{equation}
\tilde{\bm{e}}_{f,c} \coloneqq {\mathbb{D}}|_{(c,f)} (\bm b^{*}_{f} -
\bm b_{c}),
\label{dual_edge}
\end{equation}
where $\bm b^{*}_{f}$ is the barycenter of face $f$ and
$\bm b_{c} \in \mathbb R^{3}$ is an arbitrary point. Let
$\{\mathbb{R}^{\mathcal F}_{c}\}_{c \in C}$ be the family of averages of
local \emph{admissible} face reconstruction operators. Then,
%
\begin{equation}
\abs{c}\, \mathrm{col}_{(f)} {\mathbb{R}^{\mathcal F}_{c}} =
\tilde{\bm{e}}_{f,c}, \qquad \forall c \in C, \forall f \in F(c).
\label{identification}
\end{equation}
Moreover, if $K$ is polyhedral, then
$\{\mathbb{R}^{\mathcal F}_{c}\}_{c \in C}$ are also $P_{0}$-consistent
reconstruction operators.
\begin{proof}
Formula {(\ref{identification})} is a direct consequence of Equation (35)
in Proposition 3.3 of \cite{Brezzi2014MimeticSP}, the definition of averages
of local reconstruction operators in {(\ref{average})} and the definition
{(\ref{dual_edge})}. The fact that formula {(\ref{identification})} defines
$P_{0}$-consistent reconstruction operators follows by applying Lemma 3
in \cite{Pitassi2021TheRO} and the reasoning to deduce {(\ref{eqq2})} in the
proof of {Lemma~\ref{node_lemma}}.
\end{proof}
\end{proposition}
If $K$ is a curved grid then the admissible reconstruction operators in
{(\ref{identification})} do not satisfy the accuracy property
\textbf{(P1)}, and a simple example is given by the curved cube in
{Fig.~\ref{curved_cube}}. This inconsistency is the theoretical reason for the
lack of convergence of the MFD method on curved grids.

\section{Numerical results}
\label{numeric}

In this section we present numerical experiments to test the consistency
and convergence of the curved MFD method. We consider a stationary conduction
problem as a prototype example of an elliptic boundary value problem. The
stationary conduction problem in a domain $\Omega $ reads
%
\begin{subequations}
\label{problem}
%
\begin{align}
&\curl \bm E = \bm 0,
\label{eq42a}\\
&\div \bm J = 0,
\label{eq42b}\\
&\bm J = \sigma \bm E,
\label{eq42c}
\end{align}
\end{subequations}
where $\bm E$ is the electric field, $\bm J$ is the current density, and
$\sigma $ is the material parameter electric conductivity.

The boundary $\partial \Omega $ is partitioned as
$\partial \Omega = \partial \Omega ^{c} \cup \partial \Omega ^{nc}$, where
$\partial \Omega ^{c}$ is the subset of $\partial \Omega $ where electrodes
are physically placed. We assume that $\partial \Omega ^{c}$ decomposes
as the union of $N+1$ disjoint regions
$\{\partial \Omega ^{c}_{i}\}_{i=0}^{N}$, each representing a single electrode.
For each electrode $\{\partial \Omega ^{c}_{i}\}_{i=1}^{N}$ we enforce
a corresponding electromotive force $\{U_{i}\}_{i=1}^{N}$ with respect
to a reference electrode $\partial \Omega ^{c}_{0}$. On the remaining part
of the boundary $\partial \Omega ^{nc}$, we impose a homogeneous Neumann
boundary condition.

We derive a mixed variational formulation of problem {(\ref{problem})} as
follows. First, the curl-free condition $\curl \bm E = \bm 0$ allows us
to introduce the scalar potential $U$ such that $\bm E= -\grad U$.
Next,
we introduce the Sobolev space
$H(\mathrm{div},\Omega ) = \{ \bm v \in (L^{2}(\Omega ))^{3} \mid
\div \bm v \in L^{2}(\Omega )\}$, and the Sobolev spaces
$H_{0,\partial \Omega^{nc}}(\mathrm{div},\Omega ) = \{ \bm v \in H(\mathrm{div},\Omega )
\mid \bm v \cdot \bm n_{\partial \Omega} = 0 \text{ on } \partial
\Omega^{nc} \}$,
$H^{1}_{0}(\Omega )=\{ q \in L^{2}(\Omega ) \mid \grad q \in (L^{2}(
\Omega ))^{3},\, q = 0 \text{ on } \partial \Omega \}$. The variational
formulation reads as: find $\bm J \in H_{0,\partial \Omega^{nc}}(\mathrm{div},\Omega )$ and
$U \in H^{1}_{0}(\Omega )$ such that
%
\begin{subequations}
\label{eq43}
%
\begin{align}
&\int _{\Omega }\sigma ^{-1} \bm J \cdot \bm v \, dV - \int _{\Omega }U
\cdot \div \bm v \, dV = \int _{\Omega }\sigma ^{-1} \bm J_{s} \cdot
\bm v \, dV, \, \forall \bm v \in H_{0,\partial \Omega^{nc}}(\mathrm{div},\Omega ),
\label{eq1}\\
&\int _{\Omega }\bm J \cdot \grad q\, dV = 0, \, \forall q \in H^{1}_{0},
\label{eq2}
\end{align}
\end{subequations}
where the vector field $\bm J_{s} \in H_{0,\partial \Omega^{nc}}(\mathrm{div},\Omega )$ is introduced
to take into account Dirichlet boundary conditions on the electrodes surface
$\partial \Omega ^{c}$.

We now introduce the curved MFD discretization of equations
{(\ref{eq1})}, {(\ref{eq2})}. The train of thought is exactly the same of the
standard MFD method except that we employ $P_{0}$-consistent reconstruction
operators in place of the standard MFD ones. Therefore, the current density
$\bm J$ is approximated in the vector subspace of face DoFs
$\mathcal F_{0,\partial \Omega^{nc}} \subset \mathcal F$ that collects arrays $\mathbf J$ such
that ${\mathbf J}|_{(f)}=0$ for each $f \in \partial \Omega^{nc}$, and the
electric potential $U$ is approximated in the vector space of elements
DoFs $\mathcal C= \mathbb R^{\abs{C}}$.
The discrete counterpart of the divergence operator is given by the element-face
incidence matrix $\mathbb{D}$.

In order to define the inner product
$[ {\cdot } \, , {\cdot } ]^{\mathcal F}$ on $\mathcal F$, we construct
a special family
$\{ \mathbb{R}^{\mathcal F}_{\bm b_{c}} \}_{c \in C}$ of $P_{0}$-consistent
reconstruction operators and we proceed in six elementary steps as follows:
\begin{enumerate}
\item We compute the points $\{\bm b_{f}\}_{f \in F}$, defining rows of
matrix $\mathbb{R}^{\mathcal F}$ as in {(\ref{generalized_bary})}, as the
optimal solution of the following optimization problem:
%
\begin{equation}
\begin{aligned}
& \min _{f \in F} \max _{f \in F} & &
\norm{\bm b'_{f} - \bm b^{*}_{f}}_{\infty }
\\
& \text{such that} & & \{\bm b'_{f}\}_{f \in F}
\text{ are a solution of } \text{{(\ref{main_linear})}},
\end{aligned}
\label{opt_problem}
\end{equation}
where points $\{\bm b_{f}^{*}\}_{f \in F}$ are the barycenters of faces
and $\norm{\cdot}_{\infty}$ is the sup-norm.
\footnote{The minimization problem {(\ref{opt_problem})} can be rewritten as
a constrained linear optimization problem using standard manipulations
that we omit.}
%
\item We construct the family
$\{\mathbb{R}^{\mathcal F}_{\bm 0_{c}}\}_{c \in C }$ of $P_{0}$-consistent face reconstruction operators using
{(\ref{formula_op})} with points $\{\bm b_{f}\}_{f \in F}$.
\item We compute the points $\{\bm b_{c}\}_{c \in C}$ as
%
\begin{equation}
\bm b_{c} \coloneqq \frac{1}{\abs{F(c)}} \sum _{f \in F(c)} \bm b_{f},
\qquad \forall c \in C.
\label{bary_formula_algo}
\end{equation}
\item We construct the family
$\{ \mathbb{R}^{\mathcal F}_{\bm b_{c}} \}_{c \in C}$ of 
$P_{0}$-consistent reconstruction operators using
{(\ref{virtual_dual_node})} with points $\{\bm b_{c}\}_{c \in C}$ and reconstruction operators $\{\mathbb{R}^{\mathcal F}_{\bm 0_{c}}\}_{c \in C }$.
\item We construct the local mass matrices
$\{\mathbb{M}^{\mathcal F}_{c}\}_{c \in C}$ using {(\ref{local_mass})} with reconstruction operators
$\{ \mathbb{R}^{\mathcal F}_{\bm b_{c}} \}_{c \in C}$.
\item Finally, the mass matrix $\mathbb{M}^{\mathcal F}$ is constructed
by using {(\ref{global_mass})}.
\end{enumerate}
%

%
\begin{remark}
\label{rem5}
The rationale behind the definitions of points
$\{\bm b_{f}\}_{f \in F}$ in {(\ref{opt_problem})} and 
$\{\bm b_{c}\}_{c \in C}$ in {(\ref{bary_formula_algo})} stems from the geometric
interpretation of $P_{0}$-consistent reconstruction operators outlined
in {Remark~\ref{geometry_consistency}} and ensures that the dual grid structure
associated with $P_{0}$-consistent reconstruction operators
$\{ \mathbb{R}^{\mathcal F}_{\bm b_{c}} \}_{c \in C}$ is nicely staggered
with respect to the given primal grid $K$; see {Fig.~\ref{cube_patch}} and
{Fig.~\ref{sphere}} for illustrations of such nicely staggered dual grids. In
practice, this property is very important because is directly related with
optimal properties of local reconstruction operators and thus of local
inner products; see \cite[Section 5.2]{Pitassi2021TheRO}.
\end{remark}

\begin{remark}
The particular choice of the sup-norm in the objective function in
{(\ref{opt_problem})} is motivated by two facts: (i) the unit ball of $\mathbb R^{3}$ in the sup-norm is a cube; (ii) in the numerical test we will focus only on cubic
grids with randomly perturbed nodes.
\end{remark}

The mimetic weak formulation of {(\ref{eq1})}, {(\ref{eq2})} reads: find
$\mathbf J \in \mathcal F_{0,\partial \Omega^{nc}}$ and $\mathbf U \in \mathcal C$ such that
%
\begin{subequations}
\label{eq47}
%
\begin{align}
& [ {\mathbf J} \, , {\mathbf v} ]^{\mathcal F}- [ {\mathbf U} \, , {
\mathbb{D} \,\mathbf v} ]^{\mathcal C}= [ {\mathbf J_{s}} \, , {
\mathbf v} ]^{\mathcal F}, \, \forall \mathbf v \in \mathcal F_{0,\partial \Omega^{nc}}
\label{eq21}
,\\
&[ {\mathbb{D} \,\mathbf J} \, , {\mathbf q} ]^{\mathcal C}= [ {
\mathbf 0} \, , {\mathbf q} ]^{\mathcal C}, \, \forall \mathbf q \in
\mathcal C,
\label{eq22}
\end{align}
\end{subequations}
where $\mathbf J_{s} \coloneqq P^{\mathcal F}(\bm J_{s})$ is the projection
onto $\mathcal F_{0,\partial \Omega^{nc}}$ of the vector field $\bm J_{s}$.

 In what follows, we will need the following definitions. Let
$K_{h}=(N_{h}, E_{h}, F_{h}, C_{h})$ be a cubic grid with edge length
$h$ on a domain $\Omega $, and let $K=(N, E, F, C)$ be a curved grid in
$\mathfrak M(K_{h})$ as defined in Section~\ref{sufficient_condition}. Thus, each
face $f \in F$ is a piecewise linear surface made of two triangles
$t_{f}^{(1)}, t_{f}^{(2)}$. We say that a curved face $f\in F$ is
\emph{treated as a planar face} if the two triangles
$t_{f}^{(1)}, t_{f}^{(2)}$ decomposing it are considered as individual
faces. For a given $l \in \{0, \dots , \abs{F}\}$, we denote by
$K_{l}=(N_{l}, E_{l}, F_{l}, C_{l})$ the curved grid obtained by treating
as planar faces a subset of $l$ faces of $F$. Note that $K = K_{0}$, and
$K_{poly} \coloneqq K_{\abs{F}}$ is a polyhedral grid with twice the number
of faces of $K$. Finally, by saying that we use the standard MFD or DGA
method on the curved grid $K_{l}$ we mean that we use reconstruction operators
{(\ref{identification})} in Section~\ref{comparison} for each face in $F_{l}$, i.e.,
we use barycenters of faces in $F_{l}$ to define {(\ref{dual_edge})}.

\subsection{Patch test}
\label{patch}

We first test the consistency of the curved MFD method by solving a stationary
conduction problem whose solution is uniform in a domain $\Omega $.

We consider the cubic domain $\Omega = [0,1]^{3}$, and we construct a curved
grid $K \in \mathfrak M(K_{h})$ starting from an initial cubic grid
$K_{h}$ on $\Omega $. We set boundary conditions to
generate a uniform current density $\bm J$ of amplitude $1$ and directed
downwards along the vertical axis.

Using $P_{0}$-consistent consistent reconstruction operators of the curved
MFD method, we exactly reconstruct (up to machine precision) the uniform
current density $\bm J$ for each curved element of $K$. In
{Fig.~\ref{cube_patch}} we picture the resulting reconstructed vector field together
with the dual grid associated with the computed $P_{0}$-consistent reconstruction
operators on $K$. Similarly, using the approach in
\cite{Brezzi2006CONVERGENCEOM} we also reconstruct exactly the uniform
current density, but we employ three times more face DoFs.

Using the standard MFD or DGA methods on the curved grid $K$ we potentially
employ one DoF for each curved face. However, we do not achieve consistency
since the reconstructed field on each element of $K$ is not uniform. To
quantitatively measure this consistency error, we compute the relative
error on the numerical approximation of the dissipated power of the cubic
resistor $\Omega $, and for the considered curved grid $K$ we obtain a
$1 \%$ relative error.


\begin{figure}
\centering
\begin{subfigure}{.5\textwidth}
  \centering
  \includegraphics[width=0.7\linewidth]{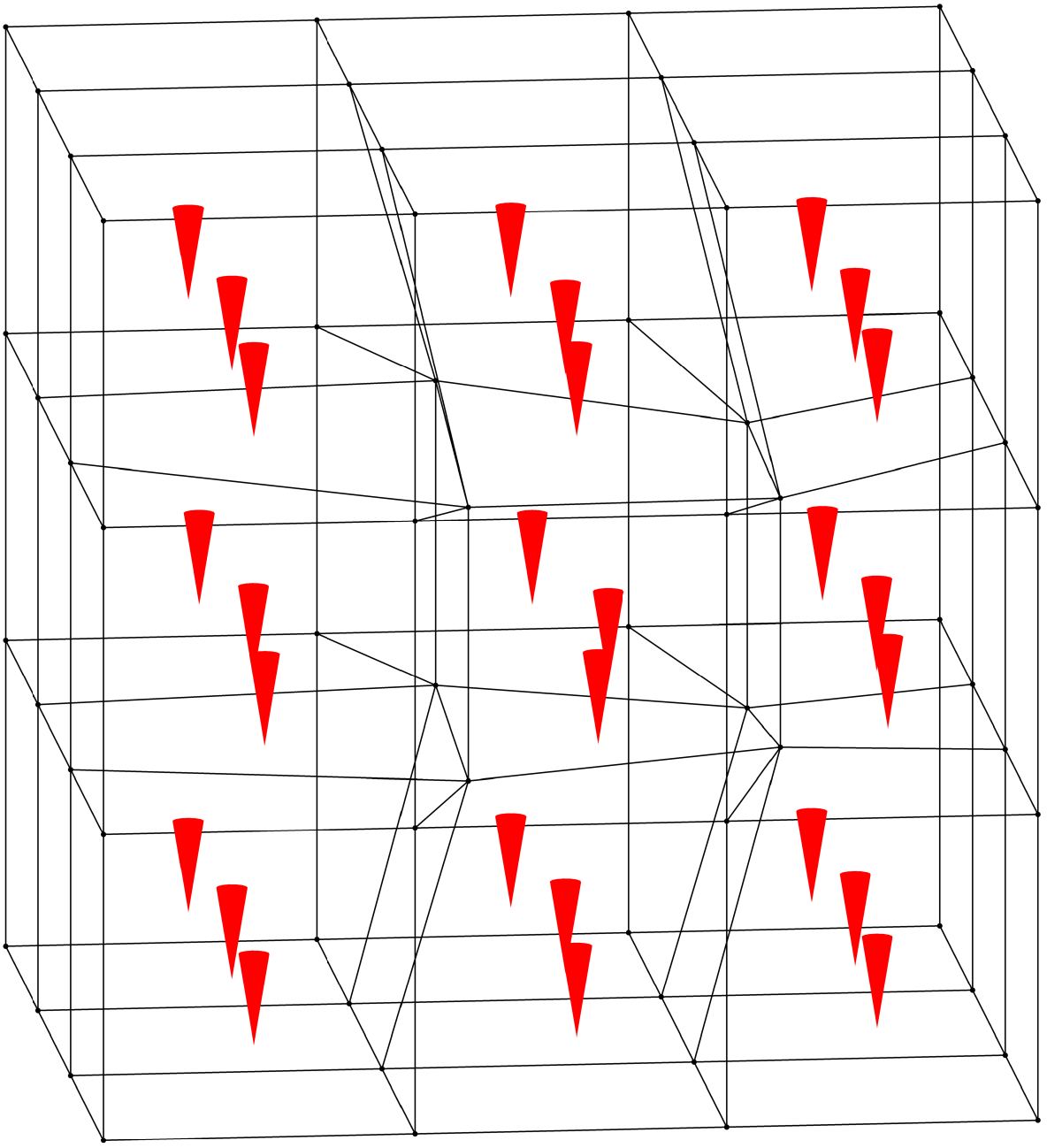}
  \label{fig:sub1}
\end{subfigure}%
\begin{subfigure}{.5\textwidth}
  \centering
  \includegraphics[width=0.7\linewidth]{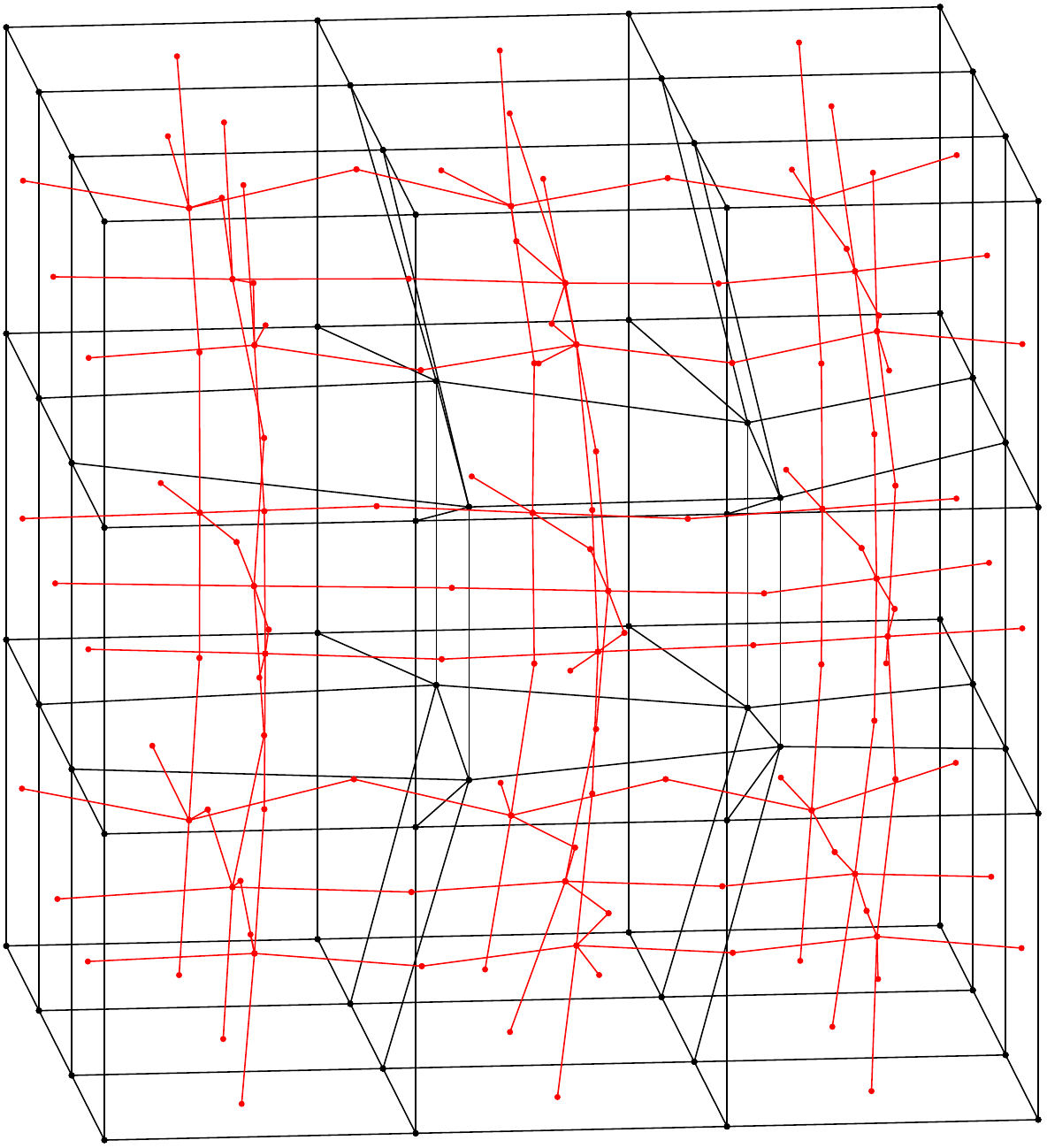}
  \label{fig:sub2}
\end{subfigure}
\caption{A curved grid partitioning the cubic resistor where all internal faces
are curved. On the left, the uniform current density reconstructed inside every
curved element; it coincides with the analytical value. On the right, in red,
the dual grid structure associated with $P_{0}$-consistent face reconstruction
operators solution of {(\ref{main_linear})}.}
\label{cube_patch}
\end{figure}

\subsection{Convergence test on randomly perturbed grids}
\label{convergence_sec}

We now test the convergence of the curved MFD method by solving a stationary
conduction problem with sufficient smooth solution on a sequence of refined curved
grids.

We consider again the cubic domain $\Omega = [0,1]^{3}$ and a sequence
of refined curved grids $K \in \mathfrak M(K_{h})$ constructed from a corresponding
sequence of cubic grids $K_{h}$ on $\Omega $ with decreasing edge length
$h$.

We measure the convergence of the discrete current density
$\mathbf J$ solution of {(\ref{eq21})}, {(\ref{eq22})} with the natural norm
$\norm{\cdot}^{\mathcal F}\coloneqq [ {\cdot } \, , {\cdot } ]^{
\mathcal F}$ induced by the discrete inner products on the space of face
DoFs $\mathcal F$. Specifically, we consider the relative error
$e^{\mathcal F}(\mathbf J)$ defined by
%
\begin{equation}
e^{\mathcal F}(\mathbf J) \coloneqq
\frac{\norm{\mathbf J - \mathbf J_{e}}^{\mathcal F}}{\norm{\mathbf J_{e}}^{\mathcal F}},
\label{relative_error}
\end{equation}
where $\mathbf J_{e}$ are DoFs in $\mathcal F$ of an exact solution
$\bm J_{e}$; thus
the quantity $e^{\mathcal F}(\mathbf J)$ provides an approximation of the
relative error in the $L^{2}$-norm on $\Omega $. In this test, we consider
the exact solution $\bm J_{e} = -\sigma \grad U_{e}$ computed with the
harmonic potential
%
\begin{equation}
U_{e}(x_{1},x_{2},x_{3}) = x_{1}^{2} - 2x_{2}^{2} + x_{3}^{2} \quad
\text{in } \Omega .
\label{harmonic}
\end{equation}

In {Fig.~\ref{convergence_cube}} we plot the convergence of
$e^{\mathcal F}(\mathbf J)$ for the curved MFD method on the sequence of
refined curved grids $K$. In addition, we also plot the convergence of
the standard MFD method applied to the corresponding sequence of polyhedral
grids $K_{poly}$ constructed by treating all faces of $K$ as planar.
Remarkably, we observe that the curved MFD method achieves the same convergence
rate but using half the number of face DoFs.

\begin{figure}
\includegraphics[scale=0.5]{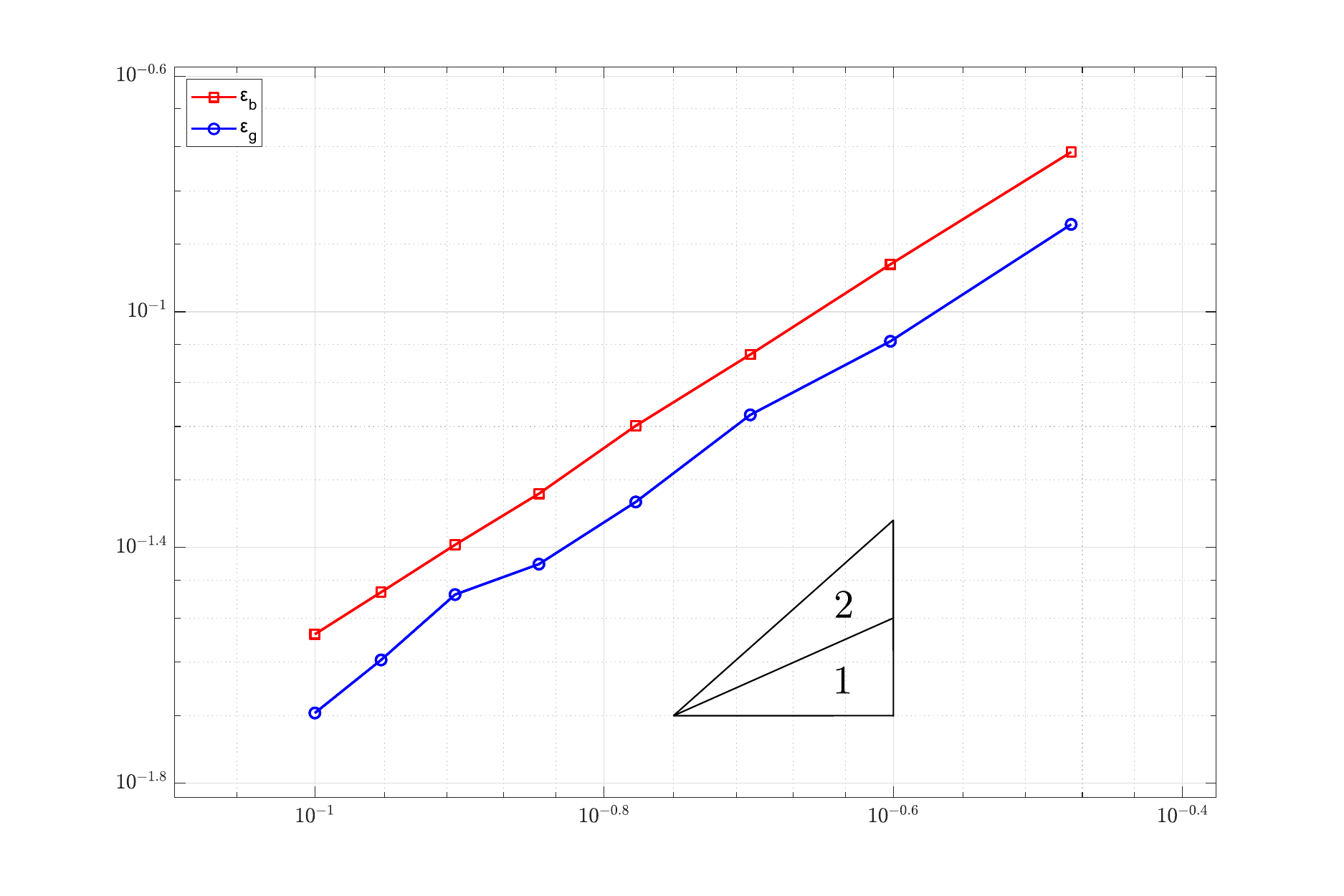}
\caption{Convergence of the relative error $e^{\mathcal F}(\mathbf J)$ in
{(\ref{relative_error})} plotted against the mesh step of the initial unperturbed
cubic grids. $\epsilon _{g}$ and $\epsilon _{b}$ represent $e^{\mathcal
F}(\mathbf J)$ for the curved MFD method and the standard MFD method where
curved faces are approximated by piecewise linear surfaces, respectively. We note that both methods achieve the same convergence rate.}
\label{convergence_cube}
\end{figure}

\subsection{Spherical resistor}
\label{sec5.3}

Next, we test the convergence of the curved MFD method on a computational domain with a curved boundary. In this way, we also test the correctness of boundary conditions
on curved boundary faces.

We consider one eight of a spherical resistor $\Omega $ illustrated in
{Fig.~\ref{curved_element}}. We enforce an electromotive force $U_{1}=1$ V between
the external curved electrode $\partial \Omega _{1}^{c}$ and the internal
one $\partial \Omega _{0}^{c}$ so that a current density $\bm J$ flows
along the radial direction.

The domain $\Omega $ is partitioned into a curved grid $K$ as schematically
illustrated in {Fig.~\ref{curved_element}}, where we note that elements of
$K$ have also curved edges. In {Fig.~\ref{convergence_sphere}} we plot the convergence
of the numerical approximation of the dissipated power versus a sequence
of refined grids of the curved MFD method. We also plot the convergence
of the standard MFD method applied to the corresponding sequence of polyhedral
grids $K_{poly}$ constructed by treating all faces of $K$ as planar;
in this case, each curved edge of $K$ is replaced with a straight edge
(segment) in $K_{poly}$.

Interestingly, for coarse grids the curved MFD method gives a better accuracy
in the approximation of dissipated power since is not affected by the geometric
error made in the approximation of the curved electrodes surface
$\partial \Omega ^{c}$ with polygons. For increasing grid sizes, this error
tends to zero and the two curves eventually coincide thus giving the same
convergence behavior.

\begin{figure}
\includegraphics[scale=0.25]{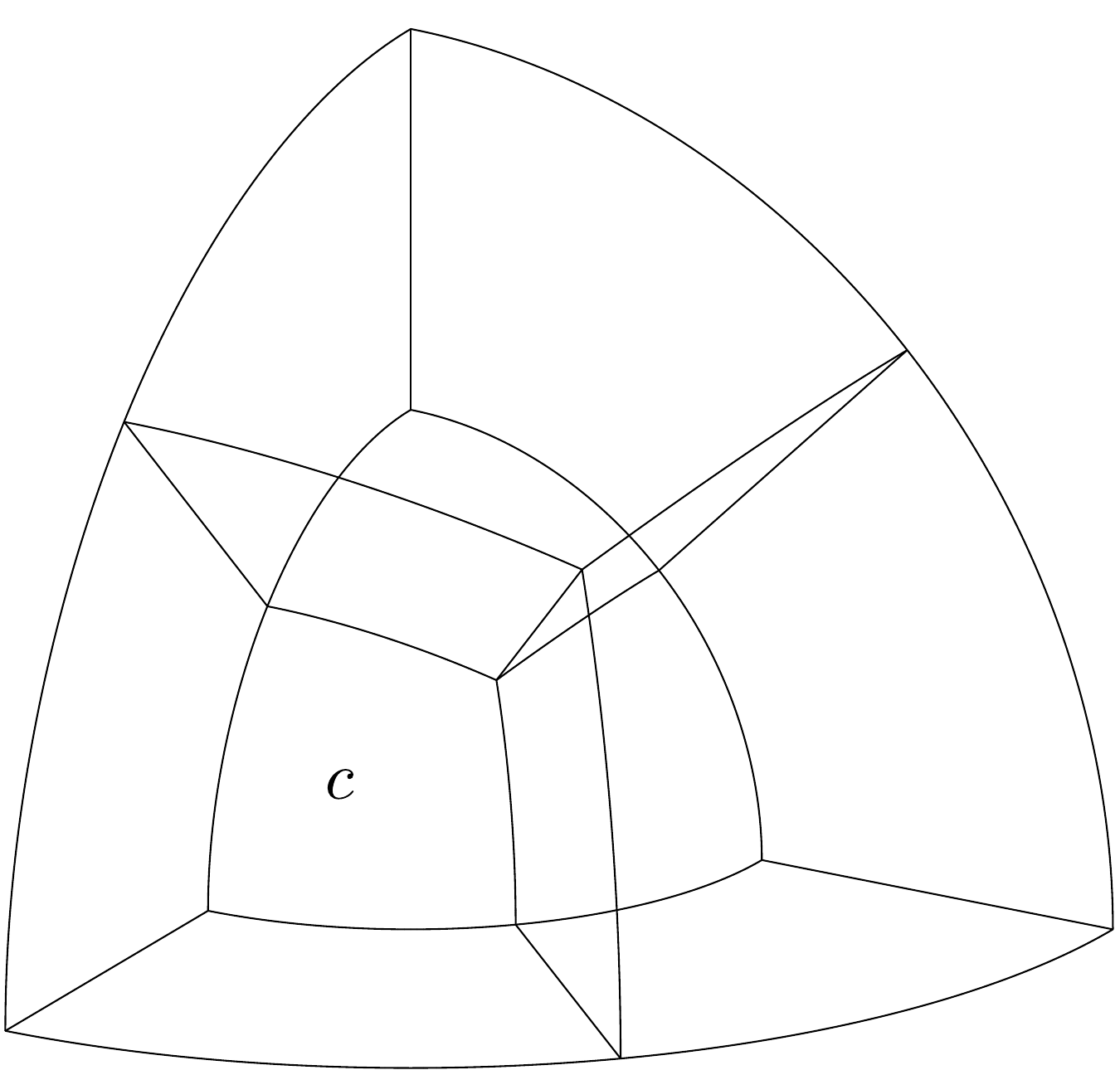}
\caption{Skeleton of a coarse curved grid partitioning one eight of a spherical
resistor. On the bottom left a curved element $c$; we note that every curved
element has both curved faces and curved edges.}
\label{curved_element}
\end{figure}

\begin{figure}
\centering
\begin{subfigure}{.5\textwidth}
  \centering
  \includegraphics[width=0.8\linewidth]{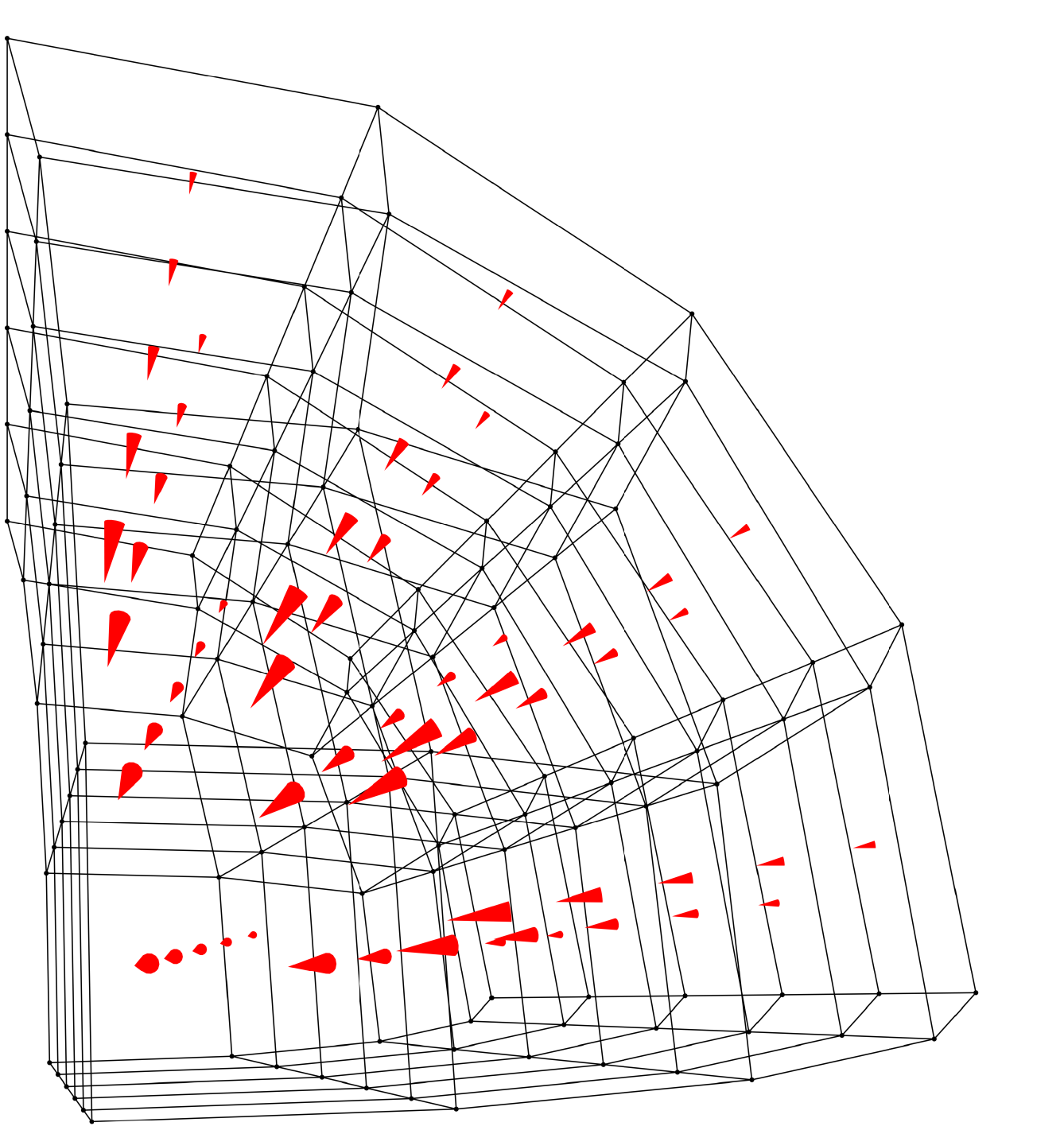}
  \label{fig:sub1}
\end{subfigure}%
\begin{subfigure}{.5\textwidth}
  \centering
  \includegraphics[width=0.8\linewidth]{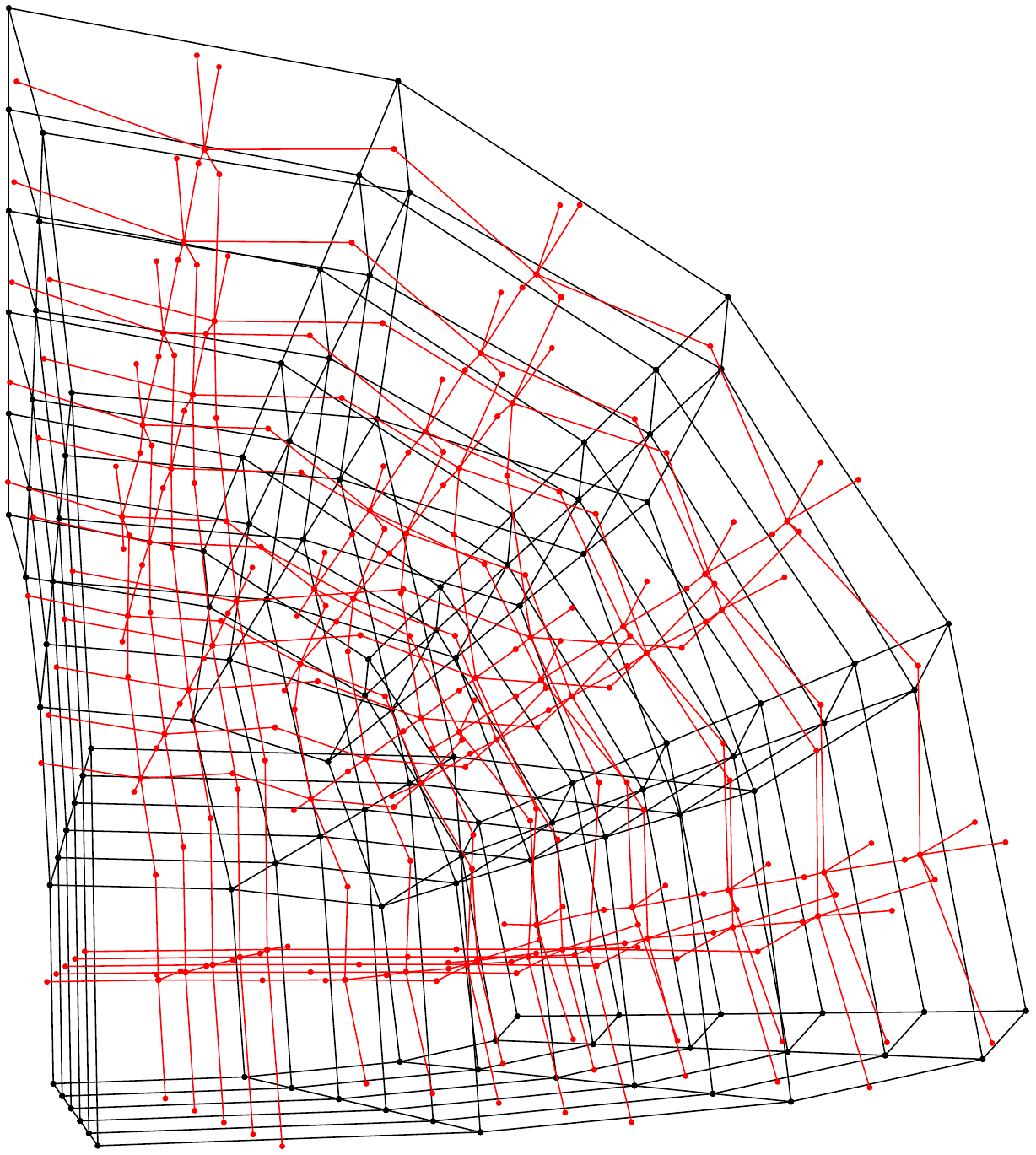}
  \label{fig:sub2}
\end{subfigure}
\caption{A curved grid partitioning one eight of a spherical resistor where all
internal faces of every element are curved just like in
{Fig.~\ref{curved_element}} (the grid nodes are connected by straight edges only
for visualization purposes). On the left, the current density reconstructed
inside every curved element. On the right, in red, the dual grid structure
associated with $P_{0}$-consistent face reconstruction operators solution of
(\ref{main_linear}).}
\label{sphere}
\end{figure}

\begin{figure}
\includegraphics[scale=0.5]{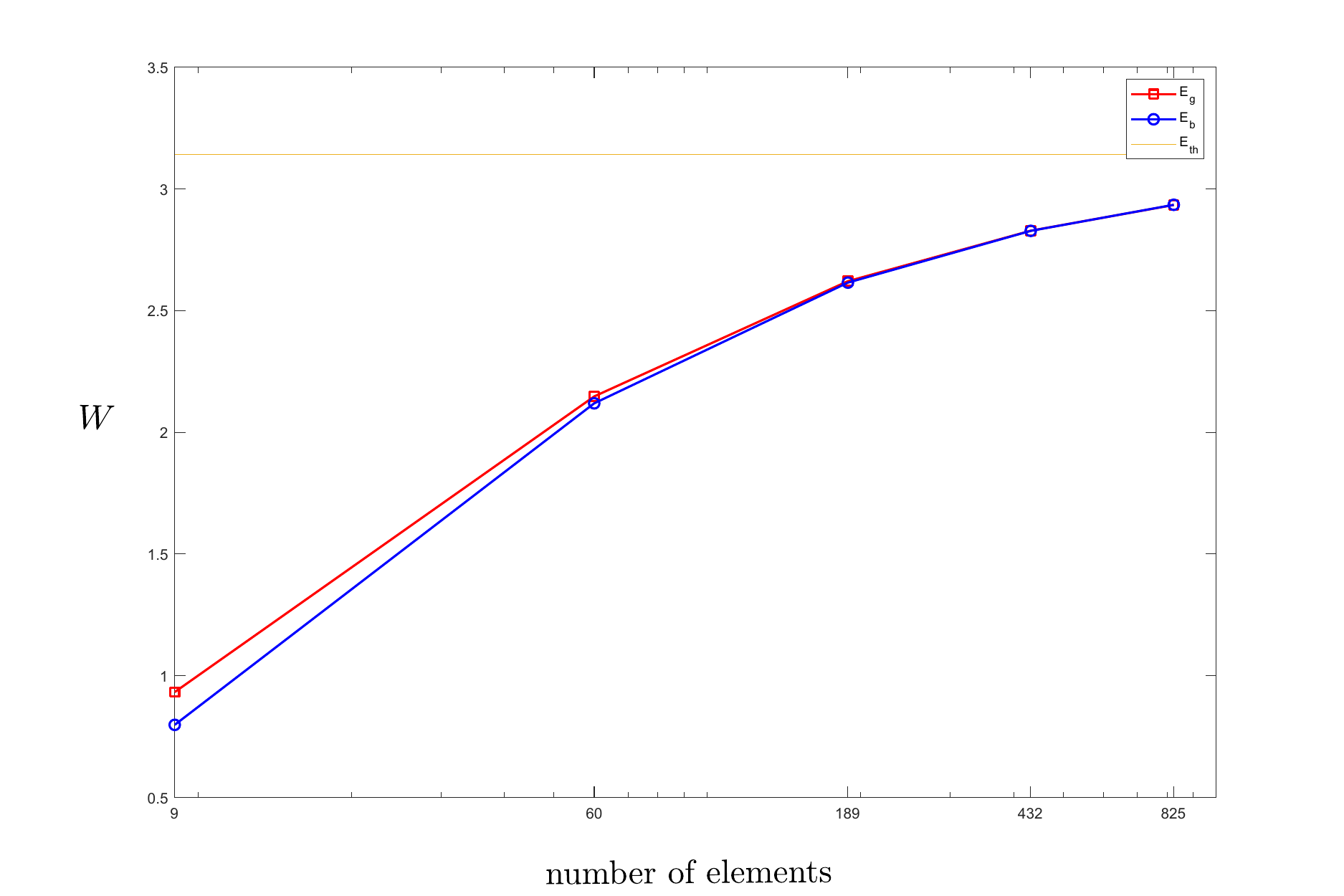}
\caption{Convergence of the numerical approximation of the dissipated power for
grids of increasing size. $E_{th}$ is the reference value equal to $\pi $ W for
the considered one eight of the spherical resistor domain $\Omega $. $E_{g}$ and
$E_{b}$ represent the numerical approximation of the dissipated power for the
curved MFD method and the standard one, respectively.}
\label{convergence_sphere}
\end{figure}

\subsection{Toroidal resistor}
\label{square_numeric_des}

Finally, we compare the curved and the standard MFD method on a same curved grid.
As observed in Section~\ref{patch}, both schemes employ
only one DoF per curved face.
However, the standard MFD method on curved grids does not converge, as due to a lack of consistency.
Nonetheless, we still want to
quantify such consistency error since may be of interest in practical applications.

We consider a toroidal resistor $\Omega$, and we construct a curved grid 
$K \in \mathfrak M(K_{h})$ obtained from an initial cubic grid
$K_{h}$ on $\Omega $.
We enforce an electromotive force of $U_{1} = 1$ V between the external lateral surface $\partial \Omega^c_1$ and the internal one $\partial \Omega^c_0$.

We compute the numerical approximation of the dissipated power of the toroidal
resistor for the curved MFD method on $K$.
In {Fig.~\ref{resistor}}, we illustrate the corresponding reconstructed
vector field together with the dual grid structure associated with the
computed $P_{0}$-consistent reconstruction operators on $K$.
We also compute the same quantity using the standard
MFD method on the curved grid $K_{l}$.

Then, we compare these two numerical values by computing their absolute difference and diving it by the analytical value of the dissipated power of the toroidal resistor.
The comparison is made by varying the number
$l$ of curved faces of $K$ that are treated as planar faces in
$K_{l}$. By varying $l$ from zero to the number of faces of $K$, we
gradually include more face DoFs.
In this way the standard MFD method on
$K_{l}$ gives a discrete problem of increasing size and accuracy. By using
the curved grid $K=K_{0}$, we obtain a $3 \%$ relative error. Instead, by
using the polyhedral grid $K_{poly}$ (i.e., treating all curved faces of
$K$ as planar), we obtain a $0 \%$ relative error.
This means that, in practice, the numerical value of the dissipated power of the standard MFD method on $K_{poly}$ is the same of the curved one on $K$.
However, we note that the number of face DoFs has
doubled.

\begin{figure}
\centering
\begin{subfigure}{.5\textwidth}
  \centering
  \includegraphics[width=1\linewidth]{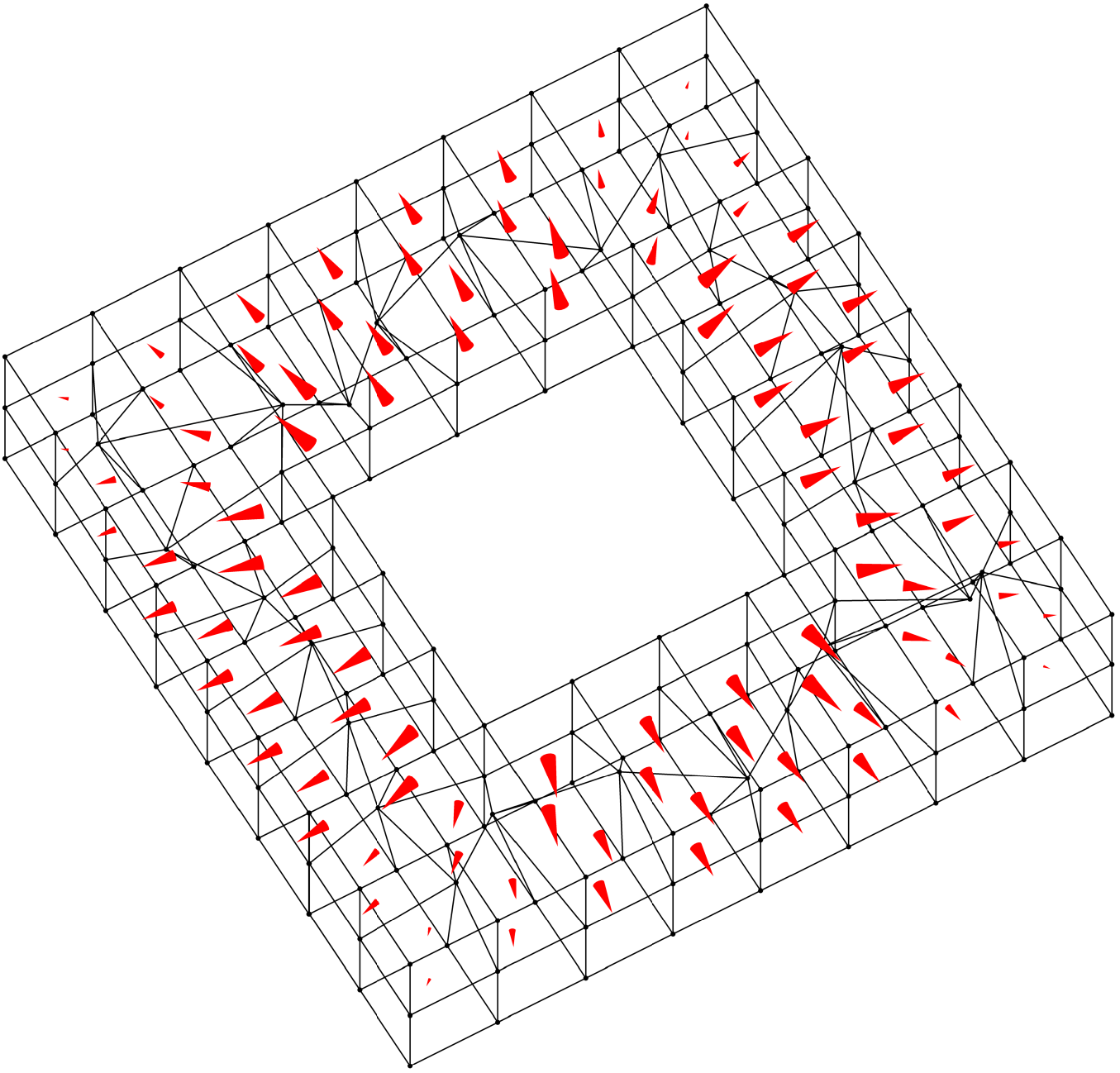}
  \label{fig:sub1}
\end{subfigure}%
\begin{subfigure}{.5\textwidth}
  \centering
  \includegraphics[width=1\linewidth]{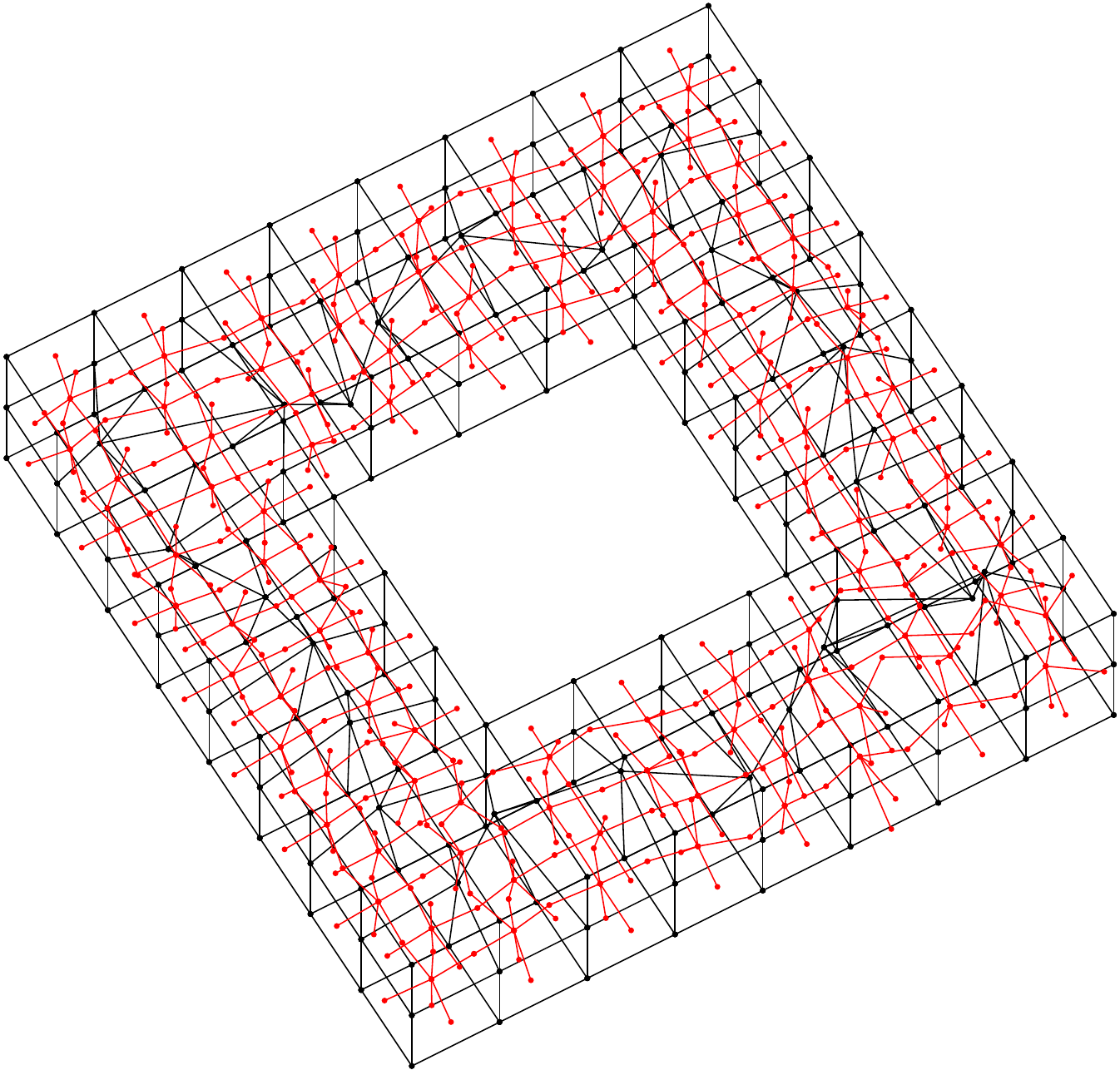}
  \label{fig:sub2}
\end{subfigure}
\caption{A curved grid partitioning the toroidal resistor where all internal faces
of each element are curved. On the left, current density reconstructed inside
each curved element. On the right, in red, dual grid structure associated with
$P_{0}$-consistent face reconstruction operators solution of
{(\ref{main_linear})}.}
\label{resistor}
\end{figure}


\section{Conclusions}
\label{conclusions}

In this work, we have proposed a new MFD method that converges on grids
with curved faces. The main novelty is that it produces a discrete problem
which is symmetric and uses only one DoF per curved face. This has been
achieved by employing the new geometrical and topological concept of
$P_{0}$-consistent reconstruction operators that generalize the standard
MFD reconstruction operators. Indeed, if the grid has no curved faces,
the proposed reconstruction operators coincide with the ones of the standard
MFD method. The consistency of the new scheme has been tested numerically,
demonstrating that the exact solution of the patch test is recovered for
discretization grids with curved faces. The convergence results show that
the curved MFD method achieves the same accuracy of the standard MFD method
where curved faces are partitioned into planar polygons but without introducing
additional DoFs. We also proposed a geometric interpretation of
$P_{0}$-consistent reconstruction operators via a dual grid structure.
As a result, the concept of dual grid has been generalized to grids with
curved faces.

A range of work is slated for future investigation, focusing on further
analysis of properties of $P_{0}$-consistent reconstruction operators.
In particular, it is our ambition to solve the characterization problem
\textbf{(Pr1)}. In this paper we have focused on problem
\textbf{(Pr2)} by proposing a possible solution strategy that, although
still needs a deep theoretical analysis, is validated from the practically
point of view by {Theorem~\ref{existence}} and by the considered numerical experiments.
Finally, we expect that the methods introduced in this paper can be extended
to high-order methods but their geometric interpretation to be more challenging.

\bibliographystyle{model1-num-names}
\bibliography{bib}

\end{document}